\documentclass[letterpaper,11pt,final]{amsart}

\usepackage[utf8]{inputenc}
\usepackage[T1]{fontenc}
\usepackage[english]{babel}
\usepackage{csquotes}

\usepackage{amssymb}
\usepackage{amsmath,amsthm}

\usepackage[notref,notcite]{showkeys}

\usepackage{cite}

\usepackage{datetime}

\usepackage[shortlabels]{enumitem}

\usepackage{todonotes}

\usepackage{yhmath}

\usepackage[pdfauthor={David Schrittesser},
            pdftitle={Definable maximal discrete sets with large continuum},
            pdfproducer={Latex with hyperref}]{hyperref}

\theoremstyle{plain}
\newtheorem{thm}{Theorem}[section]
\newtheorem{lem}[thm]{Lemma}
\newtheorem{cor}[thm]{Corollary}
\newtheorem{claim}[thm]{Claim}
\newtheorem{subclaim}[thm]{Subclaim}
\theoremstyle{definition}
\newtheorem{dfn}[thm]{Definition}
\newtheorem{exm}[thm]{Example}
\newtheorem{fct}[thm]{Fact}
\newtheorem{prop}[thm]{Proposition}
\newtheorem{rem}[thm]{Remark}
\newtheorem{rems}[thm]{Remarks}
\newtheorem{convention}[thm]{Convention}

\newtheorem*{gthm}{Galvin's Theorem}

\DeclareMathOperator{\dom}{dom}
\DeclareMathOperator{\ran}{ran}
\DeclareMathOperator{\lh}{lh}
\providecommand{\res}{\mathbin{\upharpoonright} }
\providecommand{\conc}{ \mathbin{{}^\frown}}

\DeclareMathOperator{\Power}{\mathcal{P}}

\DeclareMathOperator{\eL}{\mathbf{L}}
\DeclareMathOperator{\Vee}{\mathbf{V}}
\providecommand{\supp}{\mathrm{ supp }}
\providecommand{\setdef}{\;|\;}

\newcommand\nR{\mathrel{\ooalign{$\mathcal{R}$\cr\hidewidth$/$\hidewidth\cr}}}

\DeclareMathOperator{\otp}{otp}

\DeclareMathOperator{\proj}{proj}

\providecommand{\forces}{\Vdash}

\providecommand{\Cantorspace}{{}^\omega 2}

\DeclareMathOperator{\Sacks}{\mathbb{S}}

\DeclareMathOperator{\apo}{\mathbb{P}}

\DeclareMathOperator{\parfun}{par}
\DeclareMathOperator{\parfin}{par_{<\omega}}
\DeclareMathOperator{\init}{init}

\DeclareMathOperator{\diag}{diag}

\DeclareMathOperator{\Q}{\mathbb{Q}}
\providecommand{\Bairespace}{{}^\omega\omega}

\newcommand{\Lev}[2]{( #1 )^*_{ #2 }}
\DeclareMathOperator{\lo}{lex_0}
\DeclareMathOperator{\hi}{lex_1}

\newcommand{\Dt}[1][{}]{D^{\textbf{t}}_{#1}}

\newcommand{\Ds}[1][{}]{D^{\textbf{s}}_{#1}}

\author{David Schrittesser}
\title{Definable maximal discrete sets with large continuum}

\address{University of Toronto, %
1095 Military Trail, Toronto, ON, M1C1A4}
\email{david@logic.univie.ac.at}


\date{\today\  (\xxivtime)}


\begin{document}

\begin{abstract}
Let $\mathcal R$ be a $\Sigma^1_1$ binary relation and call a set $\mathcal R$-discrete 
\emph{iff} no two distinct of its elements are $\mathcal R$-related.
We show that in the extension of $\eL$ by iterated Sacks forcing, there is a $\Delta^1_2$ maximal 
$\mathcal R$-discrete set, 
and thus the existence of such sets is compatible with the negation of the continuum hypothesis.
As an application we find a $\Pi^1_1$ maximal orthogonal family of Borel probability measures in said extension.
The basis of this is a new Ramsey theoretic result.
\end{abstract}

\maketitle


\section{Introduction}\label{s.intro}

A map $c$ is a \emph{coloring of (pairs from) a set $X$} iff
$
c\colon [X]^2 \rightarrow \{0,1\}
$ 
where 
\[
[X]^2 = \big\{ \{ x_0, x_1 \} \setdef x_0, x_1 \in X, x_0 \neq x_1 \big\}
\]
 is the set of unordered pairs from $X$.
A set $H \subseteq X$ is \emph{homogeneous for $c$} \emph{iff} $c$ is constant on $[H]^2$. 
A major theme in  Ramsey theory is to find homogeneous sets which are large
in a given sense, such as being uncountable;
for this one needs, in general, some regularity assumption on $c$, such as being Baire measurable (with
the discrete topology on $\{0,1\}$).
Indeed,  Galvin's Theorem states that every Baire measurable coloring 
of an uncountable Polish space has a perfect homogeneous subset (see \cite[19.6,~p.~130]{kechris1995}).
Galvin's result is easily seen to be equivalent to the following statement about Sacks forcing, 
which we denote by $\Sacks$:
\begin{gthm}[in terms of Sacks forcing]\label{t.galvin.simple}
For any $p \in \Sacks$ and 
 any Baire measurable coloring $c$ of  pairs from $[p]$ 
there is 
$q \in \Sacks$ such that $q\leq p$ and $[q]$ is homogeneous for $c$.
\end{gthm}
Recall here the definition of the \emph{branch set} $[p]$ of a condition $p\in\Sacks$,
\[
[p] = \{ x \in 2^\omega \setdef (\forall n \in \omega) \; x\res n \in p \},
\]
a perfect subspaces of Cantor space $\Cantorspace$.

\medskip

Can \emph{Sacks forcing} be replaced by \emph{iterated Sacks forcing} in Galvin's Theorem?
Using notation discussed just below, we provide an answer in the following theorem; we argue in 
Fact~\ref{f.obstruction.1} 
and Example~\ref{e.obstruction.2} below that this is optimal.
\begin{thm}\label{t.galvin}
For every $\bar p$ from a dense subset of $\apo$ and every \emph{universally Baire} coloring $c$ of pairs from
$[\bar p]$  
there is $\bar q \in \apo$ such that $\bar q \leq \bar p$ and for each $\xi \in \supp(\bar q)$,
$c \res \Delta_\xi $ is continuous.
\end{thm}
Above, $\apo$ denotes a countable support iteration of Sacks forcing of arbitrary length.
In \S\ref{s.ramsey} we associate in a natural way to every $\bar p$ in a dense subset $\Dt$ of $\apo$ 
its \emph{branch space} $[\bar p]$, a perfect subspace of a countable product of copies of $\Cantorspace$.
This topologization gives meaning to the requirement that
$c$ be universally Baire.

We argue that $2^\omega$-universally Baire-ness is the optimal notion of regularity below; see the discussion after 
Fact~\ref{f.obstruction.1}.  
This is a well-behaved point-class of Baire measurable sets which includes all analytic sets and more 
(see \S\ref{s.prelim} for a review).

The role of the sets $\Delta_\xi$ is to foreclose the following obstruction to finding homogeneous sets which arises 
whenever $\lambda > 1$: For each $\xi <\lambda$, consider the coloring of pairs from ${}^{\lambda}({}^\omega 2)$ 
given by
\[
c_\xi(\{\bar x, \bar y\} ) = \begin{cases} 1 &\iff \bar x(\xi) =  \bar y(\xi),\\
0 &\text{ otherwise.}
\end{cases}
\]
The obstruction presented by this family 
leads us to partition ${}^\lambda({}^\omega 2)$:
\begin{dfn}\label{d.delta}
Given $\bar x, \bar y \in {}^\lambda({}^\omega 2)$, write $\Delta(\bar x,\bar y)$ for the least $\xi < \lambda$ such that 
$\bar x(\xi) \neq \bar y(\xi)$. Moreover, for $\xi <\lambda$ let 
\[
\Delta_\xi = \{ \{ \bar x,\bar y\}  \in  [{}^\lambda({}^\omega 2)]^2 \setdef \Delta(\bar x,\bar y)=\xi\}.
\]
\end{dfn}
\noindent
See Fact~\ref{f.obstruction.1} for further discussion.

\medskip

We apply Theorem~\ref{t.galvin} to study definable \emph{maximal discrete sets}. We use this term 
in the sense introduced by \cite{miller}:
\begin{dfn}
Let $\mathcal{R}$ be a binary relation on a set $X$. A set $A\subseteq X$ is called 
\emph{$\mathcal R$-discrete} \emph{iff} $A$ contains 
no two distinct $\mathcal R$-related elements.
By a \emph{maximal $\mathcal R$-discrete set} we mean an $\mathcal R$-discrete set which is not a proper subset 
of any $\mathcal R$-discrete set.
\end{dfn}
Discrete sets are familiar from many contexts: e.g.\ in the context of graphs (i.e.\ symmetric irreflexive relations) as 
\emph{independent sets}, or of equivalence relations, as \emph{transversals}.
Frequently one considers $\mathcal R$ arising as a \emph{compatibility relation} from a preorder: 
If $\preceq$ is a preorder, define the associated compatibility relation $\mathcal{R}_{\preceq}$ by
\[
x \mathbin{\mathcal{R}_{\preceq}} y \iff (\exists z)\;  z \preceq x \wedge z \preceq y.
\]
In such contexts, $\mathcal{R}_{\preceq}$-discrete sets are also called \emph{antichains}.

\medskip

Using the Axiom of Choice one can find maximal $\mathcal R$-discrete sets for any binary relation $\mathcal R$, 
but the existence of \emph{definable} such sets may be contentious.
In Gödel's constructible universe $\eL$ any $\Sigma^1_n$ (we use here the \emph{effective} hierarchy) 
binary relation admits a $\Delta^1_{n+1}$ maximal discrete set.
On the other hand, consider the equivalence relation $E_0$ on $\Cantorspace$, where 
$x \mathbin{E_0} y$ \emph{iff} $y(n) = x(n)$ for all but finitely many $n$.
No Baire measurable and hence no analytic or co-analytic maximal $E_0$-discrete set exists.
Moreover, if there is a Cohen real over $\eL$, no $\Sigma^1_2$ such set exists.\footnote{As another example 
let $x \mathbin{\mathcal R} y$ \emph{iff} $x \cap y$ is finite for $x$ and $y$ infinite subsets of $\omega$:
A maximal $\mathcal R$-discrete set is just a maximal almost disjoint (short: \emph{mad}) family.}

In this paper, we show that nevertheless, after forcing to add $\omega_2$-many Sacks reals to $\eL$ we can still find 
maximal discrete sets which are definable without parameters:
\begin{thm}\label{t.main}
Let $\mathcal R$ be a $\Sigma^1_1$ binary relation on an effectively presented Polish space (such as $2^\omega$ 
or $\omega^\omega$) and let $\bar s$ be generic for an iteration of Sacks forcing (of arbitrary length).
Then there is a $\Delta^1_2$ maximal $\mathcal R$-discrete set in $\eL[\bar s]$.
\end{thm}
This theorem relativizes to a parameter. 
In a previous joint work with A.\ Törnquist \cite{schrittesser-tornquist}, we proved the above for ordinary Sacks forcing 
(i.e.\ when $\lambda=1$). 
The more general Theorem~\ref{t.main} yields:
\begin{cor}
The existence of $\mathbf{\Delta}^1_2$ maximal discrete sets for all $\mathbf{\Sigma}^1_1$ binary relations is consistent with the 
negation of the continuum hypothesis (and also, simultaneously, with the negation of the Axiom of Choice).
\end{cor}

\medskip

As a more concrete application, we treat the compatibility relation associated to 
\emph{universal continuity} of  measures.
Recall that if $\mu$ and $\mu$ are (non-trivial) measures on a measurable space $X$, then one writes
 $\mu \ll \nu$ to mean that every set which is null for $\nu$ is also null for $\mu$.
Two measures $\mu$ and $\nu$ that are not compatible in $\ll$ are called \emph{orthogonal}, written $\mu \perp \nu$, 
and (maximal) $\mathcal{R}_{\ll}$-discrete sets are called \emph{(maximal) orthogonal families}. 
Orthogonal families of measures in the Polish space $P(X)$ of Borel probability measures on a Polish space $X$ 
show up in many different contexts including representation theory, ergodic theory and operator algebras.

Interest in the definability of maximal orthogonal families (or short, \emph{mof}s) originated in the following question 
posed by Mauldin: If $X$ is a perfect Polish space, is there an \emph{analytic} mof in $P(X)$?
This was answered negatively by Preiss and Rataj \cite{preiss-rataj} 
(Kechris and Sofronidis \cite{kechris-sofronidis} later gave a proof based on Hjorth's theory of turbulence).
It was shown in \cite{fischer-tornquist} that on the other hand, in $\eL$ there is a 
$\Pi^1_1$ (i.e.\ effectively co-analytic) mof in $P(X)$ whenever $X$ is an effective Polish space.

\medskip

Maximality of an orthogonal family in $P(X)$ never persists when passing to an outer model with a new real 
(an observation due to B.\ Miller, see \S\ref{s.mof}), and if there are Cohen, Random or Mathias reals over $\eL$, 
there is no $\Sigma^1_2$ 
(equivalently, no $\Pi^1_1$) mof \cite{fischer-tornquist,fischer-friedman-tornquist,schrittesser-tornquist}.
This makes it plausible that the existence of $\Pi^1_1$ mofs is essentially limited to $\eL$.
In joint work with A.~Törnquist, we showed this not to be the case \cite{schrittesser-tornquist}.
From Theorem~\ref{t.main} and the work in  \cite{schrittesser-tornquist} follows immediately the following 
strong version of this result:
\begin{thm}\label{t.mof}
The existence of a $\Pi^1_1$ mof is consistent with the negation of the continuum hypothesis.
\end{thm}

\medskip

\subsection*{Organization of the paper}
In \S\ref{s.prelim} we fix some notation regarding functions and trees and review definitions and facts regarding the 
universally Baire sets, Ramsey theory and colorings, as well as fusion for Sacks forcing and its iteration.

We prove our Ramsey theoretic result in \S\ref{s.ramsey}.
To this end, we define the set $\Dt$ of topologically determined conditions and the branch space $[\bar p]$ in 
\S\ref{s.Dt}.\footnote{\emph{Determined conditions} as well as $[\bar p]$ for  determined 
$\bar p\in \apo$ have been defined independently (and at least on the surface, differently) 
in \cite{miller-nykos}, where also a result that follows from our Corollary~\ref{c.galvin.unary} is proved. 
We thank A.\ W.\ Miller for calling this to our attention.}
This topological view of Sacks conditions is the backbone of all our arguments, and even the formulation itself
of Theorem~\ref{t.galvin}.

We then introduce the set of \emph{simple} conditions $\Ds\subseteq \Dt$, prove that $\Ds$ is dense in $\apo$ 
and at the same time \emph{continuous reading of names} for $\apo$ in \S\ref{s.continuous.reading}.
Much of this material is in some sense implicit in \cite{baumgartner-laver,groszek}; 
also compare \cite{bartoszynski,miller-nykos}. But it seemed difficult  
to prove either our Ramsey theoretic
result 
or the result about maximal discrete sets
without referring to the details of a particular form of the fusion argument (see, e.g., Definition~\ref{d.galvin.witness}) which we therefore take 
some care to work out.

We 
prove an auxiliary Ramsey theoretic result (for meager relations) in \S\ref{s.mycielski}.
Finally, we prove Theorem~\ref{t.galvin} in \S\ref{s.galvin}, in a two step argument: we show 
$\Cantorspace$-universal sets  have an auxiliary property we call $\textsf{Y}^{\apo}$-\emph{measurability}; we then show $\textsf{Y}^{\apo}$-measurable colorings can be made continuous 
by passing to a stronger condition. 

\S\ref{s.abs} presents an auxiliary result regarding absoluteness of $\mathbf{\Sigma}^1_3$-formulas for iterated
Sacks forcing, Proposition~\ref{p.elementarity}, which we are at least not aware of being available in the literature.

\S\ref{s.discrete} presents the result about maximal discrete sets, Theorem~\ref{t.main}. 
Finally \S\ref{s.mof} quickly states how to obtain Theorem~\ref{t.mof} about mofs using methods from 
\cite{fischer-tornquist} and presents the aforementioned result due to B.\ Miller (included with his permission).
We conclude with \S\ref{s.questions} listing some open questions.

\subsection*{Acknowledgments}
The author gratefully acknowledges
 the generous support from 
 Sapere Aude grant no. 10-082689/FNU from Denmark's Natural Sciences Research Council, 
the DNRF Niels Bohr Professorship of Lars Hesselholt, the Austrian Science Fund (FWF) 
through START Grant Y1012-N35, and from the Government of Canada’s New Frontiers in Research Fund (NFRF), 
through NFRFE-2018-02164.

 We also want to thank Arnie W.\ Miller and Stevo Todor{\v{c}}evi{\'c} for their comments on the present paper
 and B.\ Miller for letting us include a proof of his Theorem~\ref{t.miller}.

\section{Preliminaries}\label{s.prelim}

If $s$ is a sequence (i.e.\ a function whose domain is an ordinal), we write $\lh(s)$ for its length (said ordinal).
When $f\colon A \to B$ and $A_0\subseteq A$ we write $\{f(x) \setdef x \in A_0 \}$ as $f[A_0]$ or alternatively as 
$f''\!A_0$ to avoid confusion with the branch space.

We write ${}^A B$ for the set of functions from $A$ to $B$. 
Throughout, we shall freely identify ${}^A({}^B C)$ with ${}^{A\times B} C$ 
(in the sense that we leave it to the reader to insert the obvious identifying map when needed). 
For instance
when $\bar x \in {}^A({}^B C)$ and $B_0 \subseteq B$ 
we write $\bar x \res A \times B_0$ and treat this as an element of ${}^A({}^{B_0} C)$.
We write $\parfin(A,B)$ for the set of finite partial functions from $A$ to $B$. 
We make the same identifications for partial functions.

\medskip

Throughout this paper, let $\lambda$ be an ordinal.
We also identify 
${}^{\lambda \times \omega} 2$ and ${}^\lambda({}^\omega 2)$ as topological spaces (carrying the product topology). 
Given $s \in \parfin(\lambda\times \omega, 2)$ 
or $s \in \parfin(\lambda, \parfin(\omega, 2))$, identical in the above sense, of course
\[
N_s = \{ \bar x \in {}^{\lambda\times \omega} 2 \setdef \bar s \subseteq \bar x \}.
\]
is the basic open neighborhood in ${}^\lambda({}^\omega 2)$ defined by $s$.

\medskip

\subsection*{Codes for continuous functions}

If $X$ is a topological space with basis $\mathcal B$ and $F\in \operatorname{C}(X, {}^\omega \omega)$, 
we can consider as a \emph{code} for $F$ the function
$f \colon \mathcal B \to {}^{<\omega}\omega$ defined by letting 
\[
f(b) = \bigcup \{ s \setdef  F[b] \subseteq N_s \}.
\]
We write $f^*$ for the unique $F$ which $f$ codes, if there is such $F$.
Provided $X$ is an effective Polish space relative to a parameter $r$, being the the code of some 
$F\in \operatorname{C}(X, {}^\omega \omega)$ is arithmetical in $r$ (see \cite[2.6]{kechris1995}) 
and thus absolute for models of $\mathsf{ZF}$. 
If 
$\Vee[G]$ is a generic extension of $\Vee$ and $F \in \Vee$ is a continuous function with code $f$ we use the same 
signifier 
($f^*$ or, slightly abusing notation, $F$) for both $f^*$  as interpreted in $\Vee[G]$ and $f^*$ as interpreted in 
$\Vee$ (the restriction of the former).
In the cases arising in this article, $\{ b \in \mathcal B \setdef x \in b\}$ will be linearly ordered by 
$\subseteq$ for every $x \in X$ so that $f$ is a code \emph{iff}
\[
f^* (x) = \bigcup \{ f(b) \setdef b \in \mathcal B \wedge  x \in b \}
\]
defines a total function $f^*\colon X\to {}^{\omega}\omega$.

\medskip

\subsection*{Universally Baire sets}

Let $X$ be any Polish space. 
For a (non-meager) topological space $\Omega$,
$A\subseteq X$ is $\Omega$-universally Baire \emph{iff} for every continuous $f \colon \Omega \to X$,  
$f^{-1}[A]$ has the property of Baire (in $\Omega$);
a set $A \subseteq X$ is called universally Baire \emph{iff} $X$ is $\Omega$-universally Baire for every 
compact Hausdorff  space $\Omega$.
 The point-class of universally Baire sets is  a $\sigma$-algebra which is also closed under continuous pre-images. 
 The same holds for the $\Cantorspace$-universally Baire sets ($\Cantorspace$ denotes Cantor space).
Moreover, every \emph{absolutely $\mathbf{\Delta}^1_2$} and thus every 
$\operatorname{\sigma}(\mathbf{\Sigma}^1_1)$ set is universally Baire, and thus $\Cantorspace$-universally Baire.
See \cite{feng-magidor-woodin}, \cite[10.110,~p.~795]{woodin} and \cite{fremlin} for more details.

A map $g \colon X \to Y$ where $Y$ is a topological space is said to be $\Cantorspace$-universally Baire 
precisely if the pre-images of open sets under $g$ are $\Cantorspace$-universally Baire.

\medskip

\subsection*{Colorings}

We review some Ramsey theoretic terminology regarding colorings of pairs.
Pairs will always be from a subspace $X$ of ${}^\lambda ({}^\omega 2)$.  Let $\prec$ denote the lexicographic 
ordering of ${}^\lambda ({}^\omega 2)$---by this we mean that $\bar x \prec \bar y$ holds 
\emph{iff} for the least $\xi<\lambda$ such that $\bar x(\xi)$ differs from $\bar y(\xi)$  and for the least 
$n<\omega$ such that $\bar x(\xi)(n)$ differs from $\bar y(\xi)(n)$ we have $\bar x(\xi)(n) < \bar y(\xi)(n)$ 
(i.e.\ $\bar x(\xi)$ comes before $\bar y(\xi)$ in the lexicographic ordering of ${}^\omega 2$).
Let $\lo(\{ x_0, x_1 \})$ denote the lexicographically smaller element of the unordered pair 
$\{ x_0, x_1 \} \in [X]^2$ and  $\hi(\{ x_0, x_1 \})$ the lexicographically greater element.
As usual, we equip $[X]^2$ with the initial topology with respect to the inclusion $\iota\colon [X]^2 \to X^2$, 
$\iota(z) = ( \lo(z), \hi(z) )$
 (i.e.\ identify $[X]^2$ with the `lower half' of $X^2$).

Note that it is in some cases more natural to work with $X^2$ directly. 
In this case, we talk about symmetric maps on $X^2$ instead of colorings of pairs from $X$:
Consider a symmetric map $c^* \colon X^2 \to \{0, 1\}$. 
Any such $c^*$ clearly induces a map ${c^* \circ \iota} \colon [X]^2 \to \{0,1\}$.
Vice versa, if $c \colon [X]^2 \to \{0,1\}$, consider the symmetric map $c^* \colon X^2 \to \{0,1\}$ defined by 
$c^*(x_0, x_1) =  c( \{ x_0, x_1 \})$; let $c^*$ take value $0$ by convention on $\diag(X)$.
Baire measurability and $\Omega$-universally Baire measurability is preserved by both of these translations.
The same is not true for continuous colorings:  a continuous coloring of pairs from $X$ 
corresponds to a symmetric map which is continuous on $X^2 \setminus \diag(X)$.

\medskip

\subsection*{Trees}

All terminology regarding trees and sequences not explicitly introduced below is taken from \cite{kechris1995}.
We write $\textsc{PTrees}$ for the the $G_\delta$ subset of $\Cantorspace$ consisting of 
characteristic functions of perfect trees (as in \cite[4F, 4.32]{kechris1995}),
equipped with the subspace topology.

For a tree $T \subseteq {}^{<\omega} 2$ and $t\in {}^{<\omega} 2$,
\begin{itemize}
\item For $a \in 2$, $t \conc a$ is the element $t'$ of ${}^{<\omega} 2$ such that $t'\res \lh(t) = t$ and $t'(\lh(t))=a$.
\item $t$ is a \emph{splitting node (of $T$)} \emph{iff} there are $a,a' \in 2$ such that $a\neq a'$ and both 
$t\conc a$ and $t\conc a'$ are in $T$;
\item $T$ is \emph{perfect} \emph{iff} every node in $T$ can be extended to a splitting node of $T$;
\item Write $T_t$ for the set $\{ s \in T \setdef s \subseteq t \text{ or } t\subseteq s \}$.
\item For $n> 0$, we say $t$ is \emph{an $n$th splitting node (of $T$)} if and only $t$ is a splitting node of 
$T$ with exactly $n-1$ splitting nodes strictly below it.
That is, we start counting at $1$ (unlike Laver in \cite{laver}).
\end{itemize}

\medskip

\subsection*{Sacks forcing}

Sacks forcing $\Sacks$ is the set of perfect sub-trees of ${}^{<\omega} 2$,
ordered by $q \leq p \iff q \subseteq p$.
It admits a \emph{Fusion Lemma}, which we review below using the following terminology:
\begin{dfn}~
 \begin{itemize}
\item For $n\in \omega$ and $p,q \in \Sacks$, let us write $q \leq_n p$
to mean that $q \leq p$ and for $n>1$, 
the set of $n$th splitting nodes of $q$ is equal to 
the set of $n$th splitting nodes of $p$.
\item We say a sequence $\langle p_n \colon n\in \omega\rangle$ of Sacks conditions is a  
\emph{fusion sequence} \emph{iff} for any $m \in \omega$ there is $n_0 \in \omega$ such that for 
$(\forall n \geq n_0) \; p_n \leq_m p_{n_0}$.
\end{itemize}
\end{dfn}
\begin{fct}[Fusion for $\Sacks$]\label{l.fusion}
Any fusion sequence of Sacks conditions has a lower bound in $\Sacks$.
\end{fct}

A second important fact about Sacks forcing is that it satisfies the property in the following lemma, 
often referred to as \emph{continuous reading of names}. For a proof see e.g. \cite[3.3]{schrittesser-tornquist}.
\begin{fct}[Continuous reading of names for $\Sacks$]
Let $\dot x$ be a $\Sacks$-name for an element of ${}^{\omega}\omega$ and let $p\in\Sacks$.
Then there is $q\in\Sacks$ stronger than $p$ and a continuous function $F\colon [q]\to{}^{\omega}\omega$ such that 
$\bar q \forces_{\Sacks}  F(\bar s_{\dot G}) = \dot x $. 
\end{fct}

Our main focus will be iterations of Sacks forcing. 
Throughout, let $\lambda$ be an ordinal and let $\apo$ be an iteration of Sacks forcing with 
countable support of length $\lambda$.

We denote by $\apo_\xi$ the initial segment of length $\xi$.
Recall that $\apo$ is
the set of sequences $p\colon \lambda \rightarrow \Vee_\kappa$ 
(where $\kappa$ is some fixed, large enough ordinal) 
such that for each $\xi \in \lambda$, $p(\xi)$ is a $\apo_\xi$-name which is forced by 
$1_{\apo_\xi}$ (the trivial condition) to be a Sacks condition.
Thus if $\bar p\in\apo$, we have $\bar p \res \xi \in \apo_\xi$.
If $G$ is $\apo$-generic, $G\res\xi$ denotes $\{ \bar p \res \xi \setdef \bar p\in G\}$.

A $(\Vee, \apo)$-generic filter $G$ corresponds to an element of ${}^\lambda({}^{\omega}2)$, denoted by 
$\bar s_{G} = \langle \bar s_{G}(\xi) \colon \xi <\lambda\rangle$
such that for each $\xi<\lambda$, $\bar s_{G}(\xi)$ is a Sacks real over $\Vee[G\res\xi]$.
We emphasize again that whenever convenient we consider $\bar s_G$ to be an element of 
${}^{\lambda\times\omega}2$ via the identification of ${}^{\lambda}({}^{\omega}2)$ with 
${}^{\lambda\times\omega}2$.

\medskip

Recall that like Sacks forcing, $\apo$ admits a \emph{Fusion Lemma}:
\begin{dfn}\label{d.ifusion}
We say a sequence $\langle p_n \colon n\in \omega\rangle$ of conditions in $\apo$ is a \emph{fusion sequence} 
\emph{iff} for any $\xi \in \bigcup_{n\in\omega} \supp(p_n)$ and any  $m \in \omega$ there is 
$n_0 \in \omega$ such that for $n \geq n_0$,
$p_n \res \xi \forces$ the sets of $m$th splitting nodes of $p_n(\xi)$ and  $p_{n_0}(\xi)$ are equal.
\end{dfn}

\begin{fct}\label{l.ifusion}
Any fusion sequence of conditions in $\apo$ has a lower bound in $\apo$.
\end{fct}

When $\bar p$ is a condition in a product of Sacks forcing, $\bar p \in {}^\alpha \Sacks$, of course
\[
[\bar p] = \prod_{\xi<\alpha} [\bar p(\xi)].
\]

We sometimes decorate names in the forcing language with dots, or checks to indicate standard names, 
with the aid of helping the reader, but more often we omit such decorations.
We write $\lvert A \rvert$ for the cardinality of $A$.
Whenever $\mathbf C$ is a class or a defined object and $M$ is a model we write 
$\mathbf{C}^M$ to mean the interpretation of $\mathbf C$ in $M$.

\section{Ramsey theory of iterated Sacks forcing}\label{s.ramsey}

\subsection{Topologically determined conditions and their branch space}\label{s.Dt}
In the Ramsey theory of Sacks forcing, it is convenient to 
use the language of topology.
A simple Sacks condition $p$ corresponds to a perfect, compact set $[p] \subseteq \Cantorspace$, 
and for iterated Sacks forcing, too, it is helpful to work with analogous spaces.\footnote{Specifically,
this will help us in the proof of 
Lemma~\ref{l.baire=>Y}, leading to our version 
of Galvin's Theorem, as well as in the proof of Proposition~\ref{p.elementarity}, which 
simplifies our construction of maximal discrete sets in after iterations of Sacks forcing of length  $>\omega_1$.}

To this end the set  $\Dt\subseteq \apo$ of \emph{topologically determined} conditions is defined,  
as well the \emph{branch space} $[\bar p]$ for 
$\bar p\in \Dt$. 

The intuition is straightforward: Demand that all the Sacks conditions named by the sequence 
$\bar p\in \apo$ be described by continuous 
functions of the generic sequence of Sacks real $\bar s_G$, 
or more specifically, that $\bar p(\sigma)$ be given as a continuous function of $\bar s_G\res\sigma$. 
One additional requirement will prove to be extremely helpful\footnote{For example, in Lemma~\ref{l.Dt.hom}.}: 
namely, we require that $\bar p(\sigma)$ depend only on $\bar s_G \res (D\cap \sigma)$,
for a countable set $D \subseteq \lambda$. 

\medskip

Recall from Section~\ref{s.prelim} that we write $\textsc{PTrees}$ for the the 
effective Polish space of perfect trees (a subspace of $\Cantorspace$).

\begin{dfn}[The set $\Dt$ of topologically determined conditions]\label{d.top}~
\begin{enumerate}[label=(\alph*), ref=\alph*]
\item  We say $\bar p\in\apo$ is  
\emph{topologically determined by $H$}, abbreviated by $\bar p\in \Dt[H]$
\emph{iff}  
for some countable set
$D\subseteq \lambda$,
$H = \langle H_\sigma \setdef \sigma \in D\rangle$ is a sequence of continuous functions
\begin{equation*}
H_\sigma\colon {}^{D \cap \sigma}({}^\omega2)  \to \textsc{PTrees}
\end{equation*}
such that 
$\bar p \res \sigma \forces H_\sigma\big(\bar s_{\dot G}\res  (D\cap \sigma)\big) =  \bar p(\sigma)$
for every $\sigma \in D$, 
 and $\supp(\bar p) \subseteq D$.

\item We write $\bar p \in \Dt$ \emph{iff} there exists $H$ such that $\bar p\in \Dt[H]$ and say $\bar p$ is 
\emph{topologically determined}.
\end{enumerate}
\end{dfn}
Note that by this definition, if $\sigma = \min(D)$, $H_\sigma$ is the function with $\dom(H_0) = \{ \emptyset\}$ such that 
$H_0(\emptyset) = \bar p(\sigma)$, that is, the constant function, and 
$\bar p(\sigma)$ is forced to be a perfect tree in the ground model (for example, ${}^{<\omega}2$). 
We we prove in \S\ref{s.continuous.reading} that  $\Dt$ is dense.

\medskip

For $\bar p\in \Dt[H]$, we now define the branch space $[\bar p]^H$. Note that it depends on the choice of a 
particular sequence $H$ witnessing $\bar p \in \Dt[H]$.
\begin{dfn}[The branch space]\label{d.branchspace}~
\begin{enumerate}[label=(\alph*), ref=\alph*]
\item Given $H$ and $\bar p\in \Dt[H]$ define the \emph{branch space of $\bar p$ (with respect to $H$)} as
follows: With $D= \dom(H)$, let
\begin{equation}\label{e.branchspace.H}
[\bar p]^H = \Big\{ \bar x \in {}^D({}^\omega2) \setdef  
(\forall \sigma \in D)\; \bar x(\sigma) \in \big[H_\sigma(\bar x \res  \sigma)\big] \Big\}
\end{equation}
with the topology inherited as a subspace of ${}^{\lambda} ({}^\omega 2) $.
We write just $[\bar p]$ if $H$ can be inferred from the context.

\item\label{i.branchspace.partial} We write $[\bar p]^H_\xi$ for $\{\bar x \res \xi \res \bar x \in [\bar p]^H\}$; 
again we omit the superscript  when $H$ can be inferred from the context, or is irrelevant. 
In this case we also use the very suggestive notation
$\bar p(\xi)|_{\bar x}$ for $H_\xi (\bar x)$
where $\bar x \in [\bar p]_\xi$.
\end{enumerate}
\end{dfn}

Note that correctly interpreting $\bar x \res \emptyset$ as the empty function $\emptyset$, 
the condition on the right of \eqref{e.branchspace.H} for $\sigma=\min(D)$   just reads $\bar x(\sigma)\in [\bar p(\sigma)]$.

The space $[\bar p]^H$ (for $\bar p\in \Dt[H]$) is closed in ${}^{D\times \omega} 2$, making it 
a perfect Polish compact $0$-dimensional space; in fact, it is isomorphic to a space which is 
\emph{effectively} Polish relative to a parameter.

\medskip

Let us use the following terminology from
\cite[p.~274]{baumgartner-laver}:
\begin{dfn}
 Given $\bar p \in \apo$ and a partial function $\bar t $ from 
$\lambda \times \omega$ to $2$, we write $\bar p_{\bar t}$ for the sequence of names defined inductively 
such that for each $\xi <\lambda$,
\[
(\bar p_{\bar t})\res \xi \forces_{\apo_\xi} (\bar p_{\bar t})(\xi) = \bar p(\xi)_{\bar t}.
\]
Note that it is not necessarily the case that $\bar p_{\bar t} \in \apo$.
Therefore, we say $\bar p$ is \emph{accepts $\bar t$} precisely if $\bar p_{\bar t} \in \apo$,
or equivalently, precisely if $(\forall  \xi < \lambda)\; (\bar p_{\bar t})\res \xi \forces_{\apo_\xi} \bar t \in \bar p(\xi)$.
\end{dfn}

One can show 
$\{ N_s \cap [\bar p]^H \setdef \text{ $\bar p$ accepts $s$}\}$ 
is a topological basis for $[\bar p]^H$ and  
$\bar p \forces \bar s_{\dot G}\in\big( [\bar p]^{\check H}\big)^{\Vee[\dot G]}$.
In \S\ref{s.continuous.reading} we show that for $\bar p\in\Ds$ 
the $\bar t$ accepted by $\bar p$ essentially form a tree which we shall call $\init(\bar p)$ ($\bar t$ is accepted 
\emph{iff} it can be extended into $\init(\bar p)$). 

\medskip

\begin{rems}\label{r.change.space} 
We pause to state some rather straightforward observations regarding the relationship between the branch spaces 
of two conditions, one stronger than the other. 
Suppose $\bar p_0,\bar p_1\in\Dt$. 
\begin{enumerate}
\item 
Suppose $\bar p_i \in \Dt[H_i]$ for each $i\in \{0,1\}$.
Note $[\bar p_1]^{H_1} \subseteq [\bar p_0]^{H_0}$ \emph{iff} both 
$\dom(H_0) = \dom(H_1)$ and $\bar p_1 \leq \bar p_0$.
\item On the other hand, we are always free to enlarge the domains of $H_0$ and $H_1$ by introducing dummy variables.
Thus $\bar p_1 \leq \bar p_0$ \emph{iff}
for any $H_0, H_1$ such that $\bar p_i \in \Dt[H_i]$ for each $i\in\{0,1\}$ and such that $\dom(H_0) = \dom(H_1)$, 
we have
 $[\bar p_1]^{H_1} \subseteq [\bar p_0]^{H_0}$. 
 \item
Similarly, $\bar p_1 \leq \bar p_0$ \emph{iff} for any $H_0, H_1$ such that $\bar p_i \in \Dt[H_i]$ for each $i\in \{0,1\}$,
$\dom(H_0)\subseteq\dom(H_1)$ and the map 
$p \colon [\bar p_1]^{H_1} \to  {}^{\dom(H_0)}(\Cantorspace)$ given by $\bar x \mapsto \bar x\res \dom(H_0)$ 
is into $[\bar p_0]^{H_0}$, so that we may view it as a projection $p \colon [\bar p_1]^{H_1} \to [\bar p_0]^{H_0}$.
 \item
Suppose we have a coloring $c$ of pairs from $[\bar p_0]^{H_0}$; 
 then
$c \circ p$ is a coloring of pairs from 
$[\bar p_1]^{H_1}$.
Usually it is more convenient to choose $H_0$ and $H_1$ to have the same domain and consider the simple restriction 
$c \res [\bar p_1]^{H_1}$.
\end{enumerate}
\end{rems}

It is apparent that the same condition $\bar p \in \apo$ can give rise to different spaces $[\bar p]^H$
depending on our choice for $\dom(H)$. 
This is sometimes useful (see the next lemma) but creates ambiguity when we don't specify $H$.
We shall take care to always fix a particular $H$ when the choice matters. 
In many cases it does not matter at all, and therefore we use the following convention.

\begin{convention}
When $\bar p \in \Dt[H]$, $F\colon [\bar p]^H \to X$ is a continuous function to some space $X$,
 and $\bar x \in {}^\lambda(\Cantorspace)$, 
we write 
\[
F(\bar x)
\]
as a shorthand for
\[
F\big(\bar x \res \dom(H)\big).
\]
We shall use this convention also when $\bar p \in \Dt$ and we do not specify 
$H$ at all. 
Then too, we consider continuous functions $F\colon [\bar p] \to X$ and write
$F(\bar x)$;
by convention this should be read as $F\colon [\bar p]^H \to X$ and $F\big(\bar x\res \dom(H)\big)$ for some appropriately chosen $H$.
\end{convention}
In effect, our convention is one step short of
identifying functions $F_0$ and $F_1$ with different domains ${}^{D_0}(\Cantorspace)$ and ${}^{D_1}(\Cantorspace)$
as long as $F_i = F' \circ \proj_{D_0\cap D_1}$ for each $i\in\{0,1\}$ and some $F'$ 
with domain ${}^{D_0\cap D_1}(\Cantorspace)$, where $\proj_D(\bar x) := \bar x\res D$.
Then, we could simply always speak of 
functions $F$ defined on ${}^\lambda(\Cantorspace)$ 
which only depend on countably many components, that is, so that $F = F' \circ \proj_D$ for some countable 
$D\subseteq \lambda$.

\medskip

The next lemma is rather technical, and the reader may want to return to it after
seeing how we use it (in the proof of Claim~\ref{c.below} below,
leading to the proof Proposition~\ref{p.elementarity} and thus Theorem~\ref{t.main}, and in the proof of 
Lemma~\ref{l.baire=>Y} leading
to the proof of Theorem~\ref{t.galvin}).
I find myself unable to predict whether the reader will find the lemma surprising, trivial, or entirely opaque, 
but in any case it plays a crucial role (in the case of Lemma~\ref{l.baire=>Y}, a simpler special case suffices).

\begin{lem}\label{l.Dt.hom} Suppose $D\subseteq \lambda$ and $\pi\colon D \to \lambda$ is order preserving.
Further, suppose $\bar p \in \Dt[H]$ and $D=\dom(H)$.
\begin{enumerate}
\item\label{l.first} There are $\bar q, G$ such that $\bar q \in \Dt[G]$ and $\dom(G)=\ran(\pi)$ and a homeomorphism
\[
h\colon [\bar p]^{H} \to [\bar q]^{G}
\]
such that for any $\bar x \in [\bar p]^{H} $ and $\sigma \in D$, 
$h(\bar x)\res\pi(\sigma)$ depends only on $\bar x \res \sigma$.
\item \label{l.second}
 With $\bar q$, $G$, and $h$ as in the previous item, 
if we are furthermore given $\bar q', G'$ such that $\bar q' \in \Dt[G']$, $\dom(G') = \ran(\pi)$, 
and $\bar q' \leq \bar q$, there are $\bar p', H'$ such that $\bar p' \in \Dt[H']$ and $\dom(H') = D$
together with a homeomorphism
 \[
h'\colon [\bar p']^{H'} \to [\bar q']^{G'}
\]
which extends $h$, that is, $h' = h \res [\bar p']^{H'}$.
\end{enumerate}
\end{lem}
\begin{proof}
\eqref{l.first} 
Write $E = \ran(\pi)$. 
Define $\bar q \res \pi(\sigma)$ and 
$G\res \pi(\sigma) = \langle G_\xi \setdef \xi\in\big(E\cap\pi(\sigma)\big)\rangle$, as well as functions 
$h_\sigma$ by recursion on $\sigma$; we will do this in such a way that
\[
h_\sigma \colon  [\bar p\res\sigma]^{H\res\sigma}\to [\bar q\res\pi(\sigma)]^{G\res\pi(\sigma)}
\]
is a homeomorphism for each $\sigma$ (and they ``fit together'' to form $h$).

We start with $\sigma = \min(D)$, in which case we let $\bar q\res \sigma$ be the trivial condition in 
the iteration of Sacks forcing of length $\sigma$. 
Clearly, $G \res \pi(\sigma)$ must be the empty function;
further, must $h_\sigma$ be the unique function from $\{\emptyset\}$ to itself 
(both are notational artifacts). 

Suppose now we already have constructed $h_\sigma$, $\bar q\res\pi(\sigma)$, and $G\res\pi(\sigma)$.
For each $\bar y \in  [\bar q\res\pi(\sigma)]^{G\res\pi(\sigma)}$, let
\[
G_{\pi(\sigma)}(\bar y) = H_\sigma \circ h_\sigma^{-1}(\bar y)
\]
and define $\bar q \big(\pi(\sigma)\big)$ to be a name for 
$G_{\pi(\sigma)}\big(\bar s_{\dot G}\res \pi(\sigma)\cap E\big)$.
Further, define $\bar q(\delta)$ to be a name for the trivial condition for all $\delta$ such that 
$\pi(\sigma)<\delta<\pi(\sigma+1)$.

As a consequence, if we let $h_{\sigma+1}(\bar x) = \bar y$, where
\[
\bar y (\delta) = 
\begin{cases} h_\sigma(\bar x \res \sigma) (\delta) &\text{for $\delta \in E\cap\pi(\sigma)$,}\\
\bar x\big(\sigma) &\text{for $\delta =\pi(\sigma)$,}
\end{cases}
\]
we obtain a homeomorphism as desired.
To finish the induction, we must define 
$h_\sigma$, $\bar q\res\pi(\sigma)$, and $G\res\pi(\sigma)$ for $\sigma$ an accumulation point of $D$.
We can simply let $\bar q\res\pi(\sigma) = \bigcup_{\sigma'<\sigma}q\res\pi(\sigma')$, 
$G\res\pi(\sigma) = \bigcup_{\sigma'<\sigma}G\res\pi(\sigma')$ and let $h_\sigma$ be the function such that
$h_\sigma(\bar x) \res \pi(\sigma') = h_{\sigma'}(\bar x\res\sigma')$ for each $\sigma' \in D \cap \sigma$.

Likewise, at the end, let $h\colon [\bar p]^{H} \to [\bar q]^{G}$ be the function such that
$h(\bar x) \res \pi(\sigma) = h_{\sigma}(\bar x\res\sigma)$ for each $\sigma \in D$.

\eqref{l.second}
Since  for each $\sigma \in D$ by assumption, $h(\bar x)\res\pi(\sigma)$ only depends on $\bar x \res\sigma$, we can define 
a sequence of homeomorphisms
\[
h_\sigma \colon [\bar p]^H_\sigma \to [\bar q]^G_{\pi(\sigma)}
\]
just as in the previous proof, letting 
$h_\sigma(\bar x \res \sigma) = h(\bar x)\res\pi(\sigma)$, for each $\sigma\in D$.

Given a function $f$ and $a \notin \dom(f)$, let us momentarily write $f \conc(a,b)$ for the function extending $f$
and taking $b$ at $a$, that is, for $f\cup\big\{(a,b)\big\}$.
Consider $\sigma+1 \in D$ and let us keep $\bar x [\bar p]^H_\sigma$ fixed.
For each $x \in H_\sigma(\bar x)$, define
\[
h_{\bar x}(x)=h_{\xi+1}(\bar x\conc(\sigma,x)\} \big(\pi(\sigma)\big)
\]
---in other words, find $h_{\bar x}$ so that we may write
\[
h_{\xi+1}\big(\bar x \conc(\sigma,x)\big) = h_\sigma(\bar x) \conc\big(\pi(\sigma),h_{\bar x}(x)\big).
\]
We obtain a homeomorphism
\[
h_{\bar x} \colon [H_\sigma(\bar x)] \to [G_{\pi(\sigma)}\big(h_\sigma(\bar x)\big)]
\]
or intuitively, a continuous name for a homeomorphism from $[\bar p(\sigma)]$ to $[\bar q\big(\pi(\sigma)\big)$.
(In fact, we only require the case where as in the previous proof, $h_{\bar x}$ is the identity map;
but now that we have already discussed $h_{\bar x}$, this no longer leads to much of a simplification.)

Now let $H'_{\sigma}(\bar x)$ be the unique perfect tree $T$ such that
\[
[T] = h_{\bar x}^{-1}\big[ G'_{\pi(\sigma)}\big( h_\sigma(\bar x)\big)   \big]
\]
and let $\bar p' \in \Dt[H']$ for $H' = \langle H'_\sigma\setdef \sigma\in D\rangle$
and let $h' = h\res [\bar p']^{H'}$.
It is routine to verify that $\bar p', H'$, and $h'$ are as required.
\end{proof}

\subsection{A simpler description of iterated Sacks conditions}\label{s.continuous.reading}
We now define the set $\Ds\subseteq \Dt$ of \emph{simple conditions}
for which $[\bar p]$ will be even easier to describe;
 in Lemma~\ref{l.continuous.reading.gen} we show that $\Ds$ and hence $\Dt$ is dense.
At the same time we will establish particularly useful forms of 
``continuous reading of names'' and fusion for iterated Sacks forcing, on which the proof of
the higher dimensional version of Galvin's Theorem will be based.

\medskip

Just as a Sacks condition is a tree consisting of finite sequences from $\{0,1\}$ we define a set 
$\init(\bar p)\subseteq \parfun_{<\omega}(\lambda\times\omega,2)$ of `finite approximations'  to conditions 
$\bar p\in\apo$, which will at the same time form a basis consisting of clopen sets for the topology of 
$[\bar p]$ whenever 
$\bar p \in \Ds$ ($[\bar p]$ is well-defined as $\Ds \subseteq \Dt$).

\medskip

Let us fix the following notation:
\begin{dfn} Let $T$ be a sub-tree of ${}^{<\omega}2$.
\begin{itemize}
\item For $n\in \omega\setminus\{0\}$, write $\Lev{T}{n}$ for the set of $n$th splitting nodes and 
$\Lev{T}{\leq n}$ for $\{t'\in T \setdef (\exists t) t' \subseteq t \wedge t \in \Lev{T}{n}\}$. 
We also write $\Lev{T}{0} = \{ \emptyset\}$ for consistency. 
\item For $n\in \omega$, write $(T)_n$ for the set $t \in T$ such that $t = t' \conc a$ for some $a\in 2$ and 
$t'\in \Lev{T}{n}$ (i.e.\ the immediate $T$-successors of $n$th splitting nodes)
and write $(T)_{\leq n}$ for $\{t'\in T \setdef (\exists t) t' \subseteq t \wedge t \in (T)_n\}$.
We also write $(T)_{0} = \{ \emptyset\}$ for consistency. When there is no danger of misinterpretation, 
we leave out the brackets and write $T_n$ and $T_{\leq n}$.
\item For $n\in \omega$ and $p,q \in \Sacks$, define $q \leq_n p$
to mean that $q \leq p$ and $(q)_{\leq n} = ( p)_{\leq n}$ (or equivalently, $(q)^*_n = ( p)^*_n$).
\end{itemize}
\end{dfn}

The set $\init^n(\bar p)$, for $\bar p \in \Ds$ will generalize $p_{\leq n}$ for $p\in \Sacks$.

\begin{dfn}\label{d.init}~
Given $\bar p\in\apo$, $k,n\in \omega$ such that $k\leq n$,  and a finite or infinite sequence 
$\Sigma = \langle \sigma_l \setdef l < \alpha\rangle$ of ordinals in $\supp(\bar p)$ such that $n < \alpha \leq\omega$ and 
$\sigma_0=0$, let
\begin{gather*}\label{e.init_n}
\init^n_k(\bar p,\Sigma)  =   \{ \bar t\colon \{ \sigma_l \setdef l < n \} \cap \sigma_k \rightarrow {}^{<\omega} 2 \setdef  
(\forall l < k)\; (\bar p \res \sigma_l)_{\bar t \res \sigma_l } \forces \bar t(\sigma_l) \in \bar p(\sigma_l)_{n} \}.
\end{gather*}
We also let $\init^n(\bar p,\Sigma) = \init^n_n(\bar p, \Sigma)$ and 
$\init(\bar p, \Sigma) = \bigcup_{n\in\alpha} \init^n(\bar p,\Sigma)$.
When $\Sigma$ can be inferred from the context unambiguously, we write $\init^n_k(\bar p)$, $\init^n(\bar p)$, and
$\init(\bar p)$.
\end{dfn}
In the definition of $\init^n(\bar p)$, observe by induction on $l$ that $\bar p$ accepts $\bar t \res \sigma_l$ 
so that the definition makes sense.
Note that since we assume $\sigma_0=0$, for any $n$ and $k \leq n$ we have 
$\init^n_0(\bar p,\Sigma)=\init^0_0(\bar p,\Sigma)=\init^0(\bar p,\Sigma) = \{ \emptyset \}$. 
Also note that simply 
\[
\init^n_k(\bar p,\Sigma)  = \{ \bar t \res  \sigma_k \setdef \bar t \in \init_n(\bar p,\Sigma) \}
\]
and 
that $\init^n_k(\bar p)$ only depends on $\Sigma\res k$.

We will only be interested in $\init(\bar p)$ when $\bar p \in \Ds$; in fact, for arbitrary $\bar p\in \apo$, 
$\init(\bar p)$ may be uninteresting---for instance, it may well be that $\init(\bar p)$ contains only the trivial sequence.  
On the other hand, it will be convenient to be able to talk about 
$\init^n_k(\bar p, \langle \sigma_0, \hdots, \sigma_{n-1}\rangle)$ for arbitrary $\bar p$ 
when we prove that $\Ds$ is dense. 

\begin{dfn}[The set $\Ds$ of simple conditions]\label{d.Dcc-alt}~
\begin{enumerate}[label=(\alph*), ref=\alph*]
\item Suppose $\bar p\in\apo$. A \emph{standard enumeration of $\supp(\bar p)$} is a sequence 
\[
\Sigma = \langle \sigma_l \setdef l < \alpha\rangle
\]
where $\alpha \leq \omega$, 
$\{ \sigma_l \setdef l < \alpha \} = \supp(\bar p)$ and $\sigma_0 = 0$.

\item For $\bar p\in\apo$ we say that 
\emph{$(\Sigma, h)$ describes $\bar p$ simply}, abbreviated by $\bar p\in \Ds[\Sigma, h]$, 
\emph{iff} $\Sigma = \langle \sigma_l \setdef l < \alpha\rangle$ is a standard enumeration of 
$\supp(\bar p)$ and $h$ is a  sequence  of functions, $h=\langle h^n_k\setdef k< n <\alpha\rangle$
such that 
for any $k < n < \alpha$ it holds that $\init^n_k(\bar p, \Sigma) = \dom(h^n_k)$ and
for any $\bar t \in \init^n_k(\bar p, \Sigma)$ 
\[
 (\bar p \res \sigma_{k})_{\bar t} \forces_{\apo_{\sigma_k}} h^n_k(\bar t) =  \bar  p(\sigma_{k})_{\leq n}.
\]
\item We say that 
\emph{$\bar p$ is simple with respect to $\Sigma$}, abbreviated by $\bar p\in \Ds[\Sigma]$, \emph{iff}
there exists a sequence $h$ such that $\bar p \in \Ds[\Sigma, h]$. 
We say that \emph{$\bar p$ is simple}, abbreviated by $\bar p\in \Ds$, if and only if
there exists $\Sigma$ such that $\bar p \in \Ds[\Sigma]$.
\end{enumerate}
\end{dfn}

We pause to note that 
when $k=0$, as $\apo_{\sigma_0}$ is just the trivial forcing and $\init^n_0(\bar p) = \{ \emptyset \}$, trivially
$h^n_0$ is the function with domain $\{\emptyset \}$ taking the value $\bar p(0)_{\leq n}$. 
Further note that for $\bar p \in \Ds[\Sigma]$, there is precisely one $h$ such that $\bar p \in \Ds[\Sigma, h]$.

\medskip

Notation will be greatly simplified by the following convention:
\begin{convention}\label{convention}
Suppose $\bar p\in \Ds[\Sigma,h]$ where $h=\langle h^n_k\setdef k< n <\alpha\rangle$.
Whenever $\bar a$ is a partial function from $\lambda\times\omega$ to $2$ and $k < n < \alpha$ 
and there is some $\bar s \in \init^n_k(\bar p)$ 
such that $\bar s \subseteq \bar a$, 
we define $h^n_k(\bar a)$ to be $h^n_k(\bar s)$ (this is well-defined 
as there can clearly be at most one such $\bar s$).
\end{convention}

For the readers convenience we prove the following rather obvious fact:
\begin{lem}\label{d.pspace} It holds that $\Ds \subseteq \Dt$.
\end{lem}
\begin{proof}
Given $\Sigma= \langle \sigma_l \setdef l < \alpha\rangle$, 
$h=\langle h^n_k\setdef k< n <\alpha\rangle$ and $\bar p \in \Ds[\Sigma, h]$, we define partial functions
\[
H_{\sigma_k} \colon {}^{\sigma_k} ({}^\omega 2) \to \textsc{PTrees}
\]
for each $k < \alpha$ by
\begin{multline*}
  H_{\sigma_k}(\bar x) = 
  \bigcup \{ h^n_k (\bar x) \setdef k < n < \lh(\Sigma)\wedge\bar t \in \init^n_k(\bar p)\wedge\bar t \subseteq \bar x \} =\\
   \bigcup \{ h^n_k (\bar x) \setdef k < n < \lh(\Sigma)  \}.
\end{multline*}
Note that ${}^{\sigma_0} ({}^\omega 2) = \{\emptyset \}$; the second line makes use of Convention~\ref{convention}.
Let $H= \langle H_\xi \setdef \xi \in \supp(\bar p)\rangle$;
clearly, $H_{\sigma_k}$ is defined on $[\bar p]_{\sigma_k}^{H}$ as in Definition 
\ref{d.branchspace}\eqref{i.branchspace.partial} 
and is just 
$(\bigcup_n h^n_k )^*$ in the sense of \S\ref{s.prelim}. 
Finally, 
 $\bar p \in \Dt[H]$.
\end{proof}
This justifies introducing the following notation:
\begin{dfn}
If $\bar p \in \Ds[\Sigma,h]$ write $[\bar p]^\Sigma$ for $[\bar p]^H$, where 
$H$ is obtained as in the previous proof (noting again that $h$ is uniquely determined).
We will often drop the superscript and simply write $[\bar p]$ since $\Sigma$ will be clear from the context.
\end{dfn}
Note that if $\bar p \in \Ds[\Sigma,h]$ we also have 
\begin{multline*}
[\bar p]^\Sigma  = \{  \bar x \in {}^{\supp(\bar p)\times \omega} 2 \setdef 
\text{ For any $s \in \parfin(\supp(\bar p)\times \omega, 2)$ such that $s \subseteq \bar x$,  }\\
\text{there exists $n\in\omega$ and $\bar t \in \init_n(\bar p)$ such that $\bar s \subseteq \bar t$ and
$\bar t \subseteq \bar x$ } \}.
\end{multline*}
Thus, when $\bar p\in \Ds$, there is a simple one-to-one correspondence 
between $[\bar p]$ and the set of branches through 
$\bigcup_{n\in\omega} \init_n(\bar p)$, where we consider the latter as a tree w.r.t.\ the ordering given by $\subseteq$.
Furthermore, clearly $\{ N_s \cap [\bar p] \setdef s \in \init(\bar p)\}$ is a topological basis for the space $[\bar p]$.

\medskip

We now give a version of continuous reading for names for iterated Sacks forcing; at the same time we prove 
that $\Ds$ (and thus $\Dt)$ is dense in $\apo$.
\begin{lem}\label{l.continuous.reading.gen}
Let $\dot x$ be a $\apo$-name for an element of ${}^\omega \omega$ and $\bar p \in \apo$.
Then we can find $\bar q \in \Ds$ such that $\bar q \leq \bar p$ together with 
a map $f\colon \init(\bar q) \to {}^{<\omega}\omega$ such that for every $n\in \omega$ and $\bar s \in \init^n(\bar q)$ 
\[
(\bar q)_{\bar s} \forces \dot x \res n = f(\bar s).
\]
In particular, 
$f$ codes a continuous map
$
f^* \colon [\bar q] \to {}^\omega \omega
$
(both in $\Vee$ and $\Vee[G]$) 
such that $\bar q \forces \dot x = f^*(\bar s_{\dot G})$  (see \S\ref{s.prelim}).
\end{lem}
We introduce some more terminology that will help us build \emph{fusion sequences} in the proof of 
Lemma~\ref{l.continuous.reading.gen}:
\begin{dfn}\label{d.order}

Given conditions $\bar p, \bar q \in\apo$, $n\in\omega$ and a sequence 
$\Sigma= \langle \sigma_k \setdef k\in\alpha\rangle$ of ordinals in 
$\supp(\bar p)$ with $n \leq \alpha \leq \omega$ we write 
$\bar q\leq^{\Sigma}_n \bar p$ exactly if $\bar q\leq \bar p$ and for every $k < n$, 
\[
\bar q \res \sigma_k \forces \bar q(\sigma_k)_n = \bar p(\sigma_k)_n.
\]
We also write $\leq^{\sigma_0, \hdots, \sigma_{n-1}}_n$
for $\leq^{\Sigma}_n$ as the relation only depends on $\Sigma \res n$.
\end{dfn}
\noindent
The reader should note that $\leq^\Sigma_0$ is just $\leq$.

\medskip

Now we are ready to give the fusion argument that proves the lemma.

\begin{proof}[Proof of Lemma~\ref{l.continuous.reading.gen}]
Let $\bar p\in\apo$ and a $\apo$-name $\dot x$ for an element of ${}^\omega \omega$ be given; 
we shall find a stronger condition $\bar q$ together with sequences $\Sigma$ and $h$ witnessing 
$ \bar q\in D_{\Sigma,h}$ and a map $f$ as in  the statement of the lemma.

Let $\bar p_0 = \bar p$ and $f(\emptyset) = \emptyset$ and build a fusion sequence of conditions 
$\langle \bar p_n \setdef n\in\omega\rangle$ such that $\bar p_0 \geq \bar p_1 \geq \bar p_2 \hdots$, as follows:
Fix a standard enumeration 
$\Sigma^0=\langle\sigma^0_0, \sigma^0_1, \sigma^0_2, \hdots \rangle$ of $\supp(\bar p_0)$.
For every further step $n >0$ in the construction of the fusion sequence, after having obtained $\bar p_n$ 
we shall also fix an enumeration 
$\Sigma^n=\langle \sigma^n_k \setdef n\in\omega\rangle$ of $\supp(\bar p_{n})\setminus \supp(\bar p_{n-1})$ . 

At the end of the construction 
we will obtain a standard enumeration $\Sigma$ for $\bar q$ as follows:
Let $d \colon \omega\times \omega \to \omega$ be the well-known bijection given by:
\[
d(n,k)= \frac{(n+k+1)(n+k)}{2}+n.
\]
Then we shall let $\Sigma = \langle \sigma_l \setdef l\in\omega\rangle$ be defined by
\[
\sigma_{d(n,k)} = \sigma^n_k.
\]
In other words, $\Sigma$ will enumerate $\{ \sigma^n_k \setdef n,k\in \omega \}$ 
by the well-known \emph{diagonal counting procedure}.

Most importantly, our construction will be set up so as to guarantee that
\[
\bar p_0 \geq^{\sigma_0}_1 \bar p_1 \geq^{\sigma_0, \sigma_1}_2 
\geq \hdots \bar p_{n} \geq^{\Sigma\res n+1 }_{n+1} \bar p_{n+1}\hdots
\]
We now give the details of the successor step of the construction.
Assume $n\in\omega$ and we have already constructed $\bar p_{n'}$ and $\Sigma^{n'}$ for $n' \leq n$.
The reader should note that we have already determined at least the first $n+1$ elements 
$\langle \sigma_0, \hdots, \sigma_n \rangle$ of $\Sigma$: 
For we have determined $\sigma_{d(n',k)}$ for all $k\in \omega$ and all 
$n' \leq n$, so the first index which has not yet been assigned a value in $\Sigma$ is 
$d(n+1,0)$ and $n +1 < d(n+1,0)$ (determining in fact $\sigma_{n+1}$ but this is not relevant).

We also assume by induction that we have already defined
$h_k^{n'}$ for every pair $n',k \in \omega$ with $n' \leq n$ and $k < n'$ as well as
a $f(\bar s)$ for $\bar s \in \bigcup_{n' \leq n} \init^{n'}(\bar p_{n'})$ (these assumptions are trivial if $n=0$). 
We will now find $\bar p_{n+1}$, $h_k^{n+1}$ for $k<  n+1$ and define $f$ on $\init_{n+1}(\bar p_{n+1})$.
\begin{claim}\label{c.D-cc-succ}
Let $\bar q \in \apo$ such that $\bar q \leq^{\sigma_0, \hdots, \sigma_n}_{n+1} \bar p_n$, let 
$k\in \omega$ such that $k\leq n+1$ and let $\bar t\in \init_k^{n+1}(\bar q, \langle \sigma_0, \hdots, \sigma_k\rangle)$. 
There is
$\bar q' \in \apo$ such that $\bar q' \leq^{\sigma_0, \hdots, \sigma_n}_{n+1} \bar q$ and such
that: 
\begin{itemize}
\item If $k < n+1$, for some tree $T\subseteq {}^{<\omega}2$ we have 
\[
(\bar q \res \sigma_{k})_{\bar t} \forces \bar q(\sigma_{k})_{\leq n+1} = \check T,
\]
\item If $k=n+1$, i.e.\  $\bar t \in \init^{n+1}(\bar q)$,
for some $s \in {}^{n+1} \omega$ we have 
\[
(\bar q )_{\bar t} \forces \dot x \res \check n +1 = \check s.
\]
\end{itemize}
\end{claim}
\textit{Proof of Claim~\ref{c.D-cc-succ}}.
This is obvious: 
As $\bar t\in \init^n_k(\bar q, \langle \sigma_0, \hdots, \sigma_k\rangle)$, 
$\bar q$ accepts $\bar t$; find $\bar q^* \in \apo$ such that $\bar q^* \leq \bar q_{\bar t}$
and $\bar q^*\res\sigma_{k}$ decides $\bar q^*(\sigma_{k})$ in case $k < n+1$, i.e.\ for some $T$ we have
$\bar q^*\res\sigma_{k}\forces\bar q^*(\sigma_{k})=\check T$---or respectively, if $k = n+1$, such that
$\bar q^*$ decides $\dot x \res n+1$. 

Now let $\bar q'$ be any condition such that $\bar q' \leq^{\sigma_0, \hdots, \sigma_n}_n \bar q$
and $\bar q'_{\bar t} = \bar q^*$.
To see that such a condition $\bar q'$ indeed exists, construct $\bar q'$ by induction on $\xi<\lambda$ as follows:
Assume we have constructed $\bar q'\res \xi$ such that $\bar q'\res\xi$ accepts $\bar t\res\xi$. 
Find $\bar q'(\xi)$ so that
$(\bar q' \res \xi)_{\bar t\res\xi} \forces \bar q'(\xi) = \bar q^*(\xi)$ and
whenever $\bar r \in \apo_\xi$ is incompatible with $(\bar q' \res \xi)_{\bar t\res\xi}$,
we have $\bar r \forces \bar q'(\xi) = \bar q(\xi)$.
Clearly, $\bar q'\res\xi+1$ accepts $\bar t\res\xi+1$, verifying the induction hypothesis.
\hfill \qed{\tiny Claim~\ref{c.D-cc-succ}.}

\medskip

Applying the claim finitely many times for each $k \leq n +1 $ and each 
$\bar t\in \init^{n+1}_k(\bar p_n,\Sigma\res n)$ in turn,
we construct a finite descending chain of conditions whose last condition $\bar p_{n+1} \in \apo$ satisfies that  
$\bar p_{n+1} \leq^{\sigma_0, \hdots, \sigma_n}_{n+1} \bar p_n$
and
for any  $k< n+1 $ and 
$\bar t\in \init^k_n(\bar p_n,\Sigma\res n)$, there is $T = T(k,\bar t)$  such that
\[
(\bar p_{n+1}\res\sigma_k)_{\bar t}\forces \bar p_{n+1}(\sigma_k)_{\leq n +1} = \check T
\]
and for any 
$\bar t\in \init^n(\bar p_n,\Sigma\res n)$,
there is $s= s(\bar t)$ such that
\[
(\bar p_{n+1})_{\bar t}\forces \dot x \res\check n = \check s.
\]
Define $h_k^{n}$, with domain $\init^n_k(\bar p_n,\Sigma\res n)$
 by $h_k^{n} (\bar t) = T(k,\bar t)$ and for $\bar t \in \init^n(\bar p_n, \Sigma\res n)$, define $f(\bar t) = s(\bar t)$.

It is clear that the sequence 
$\langle \bar p_n \setdef n\in\omega\rangle$ is a fusion sequence in the sense of Definition~\ref{d.ifusion}.
By Lemma~\ref{l.ifusion}, let $\bar q$ be its greatest lower bound.
As promised, $\Sigma$ is a standard enumeration of $\supp(\bar q)$. 

For each $n \in \omega$ we have
$\bar q \leq^\Sigma_{n+1} \bar p_{n}$, i.e.\ for each $k \in \omega$ such that $k< n+1$,
\[
\bar q \res \sigma_k \forces \bar q(\sigma_k)_{\leq n+1} =  \bar p_{n+1}(\sigma_k)_{\leq n+1}.
\]
It follows that  for every $n \in \omega$,
$
\dom(h_k^{n+1})=\init^k_{n}(\bar p_n,\langle \sigma_0, \hdots, \sigma_{n-1}\rangle)= \init^k_{n}(\bar q,\Sigma).
$
 Likewise,
$
\dom(f) = \init(\bar q).
$
By construction of the functions $h_k^n$, we have $\bar q\in \Ds[\Sigma, h]$ where 
$h= \langle h_k^n \setdef  n,k\in \omega, k < n \rangle$.
\end{proof}
\begin{rem}\label{r.Ds.supp} 
The proof also shows the following: When $\bar p \in \Dt[H]$,
given a completely arbitrary standard enumeration $\Sigma$ of $\dom(H)$ we can find a stronger condition 
$\bar q \in \Ds[\Sigma]$.
We find this interesting, but we shall never need this observation.
\end{rem}

\subsection{The Polarized Mycielski's Theorem for product of Sacks forcing}\label{s.mycielski}

Before we prove Galvin's Theorem for iterated Sacks forcing, we go back to product Sacks forcing: 
The following result generalizes Mycielski's Theorem (see \cite[Theorem~1,~p.~141]{mycielski} or 
\cite[19.1,~p.~129]{kechris1995}) to infinite sequences of elements of Cantor space. 
It will be crucial in proving Galvin's Theorem for iterated Sacks forcing.
Ramsey theory for infinite sequences was also studied in \cite{todorcevic} and \cite{louveau-ea}. 
\begin{thm}\label{t.myc}
Let $B$ be a comeager subset of $X= {}^\omega({}^\omega 2)$.
Then there is a sequence $C_n$, for $n\in\omega$, of perfect subsets of ${}^\omega 2$ such that
$\big( \prod_{n\in\omega} C_n \big)\subseteq  B$.
\end{thm}
The proof is an elaboration of the argument from \cite{dorais}; a proof could also be based on 
methods from \cite{todorcevic}.
\begin{proof}
For this proof, let $\Q$ be the partial order consisting of pairs $q=( n_q, f_q)$ where $n_q \in \omega$ and
$f_q\colon {}^{\leq n_q} 2 \to {}^{<\omega} 2$ preserves $\subseteq$;
ordered by $q \leq p$ \emph{iff} $n_q \geq n_p$ and 
$f_q \supseteq f_p$.
Clearly, $\Q$ as a forcing adds a perfect subset of Cohen reals to ${}^\omega 2$.
Further, let $\bar \Q = \prod^{<\omega}_\omega \Q$, the finite support product of $\omega$ many copies of $\Q$.
We shall presently find a filter for $\bar \Q$ in the ground model, meeting countably many dense sets, 
to give us the desired sequence of perfect sets.

Let $\bar \Q^*$ denote the set of $\bar p \in \bar \Q$ such that for some $n, m\in\omega$, we have $n\leq m$ and
\begin{equation}\label{dense}
\begin{gathered}
m = \supp(\bar p),\\
(\forall k \in n) \; n_{\bar p(k)} = n, \\
(\forall k \in m \setminus n) \; n_{\bar p(k)} = 1.
\end{gathered}
\end{equation}
Clearly, $\bar \Q^*$ is dense in $\bar \Q$.

For $\bar p \in\Q^*$, write $n_{\bar p}$ and $m_{\bar p}$ for the unique $n$ and $m$ 
such that \eqref{dense} holds.
Moreover, for $\bar p \in \bar \Q^*$ and $\bar s \in {}^n({}^n 2)$, where $n=n_{\bar p}$, we shall write 
$f_{\bar p} (\bar s)$ for the function with domain $m_{\bar p}$ defined by
\[
f_{\bar p} (\bar s) (k) = \begin{cases} f_{\bar p(k)}(\bar s(k)) &  k < n,\\
f_{\bar p(k)}(\emptyset) & n \leq k < m_{\bar p}.
\end{cases}
\]
For the following two claims, let $O$ be open dense in $X$.
\begin{claim}
Suppose $\bar p \in \bar \Q^*$, $n = n_{\bar p}$ and
$\bar s\in {}^n({}^n 2)$.
Then there is $\bar q \leq \bar p$ in $\bar \Q^*$ with $n_{\bar q} = n_{\bar p}$  such that
$N_{f_{\bar q} (\bar s)} \subseteq O$.
\end{claim}
\begin{proof}
As $O$ is open dense in $X$, we may pick  an extension $\bar t' \in {}^{<\omega}({}^{<\omega} 2)$ of 
$f_{\bar p} (\bar s)$
such that $N_{\bar t'} \subseteq O$. Assume $\dom(t') = m \in \omega$.
Now let $\bar q$ be any extension of $\bar p$ in $\bar \Q^*$ such that $n_{\bar q} = n_{\bar p}$, $m_{\bar q} = m$,
 $f_{\bar q(k)}(\bar s(k))= \bar t'(k)$ for $k < n_{\bar q}$ and
$f_{\bar q(k)}(\emptyset) = \bar t'(k)$ for $n \leq k < m$.
In other words,
$f_{\bar q} (\bar s) = \bar t'$.
\end{proof}
\begin{claim}
The set
$D_O = \{ \bar q \in \bar \Q \setdef (\forall \bar s \in {}^{n_{\bar q}}({}^{n_{\bar q}} 2))\;  
N_{f_{\bar q}(\bar s)} \subseteq O \}$
is dense in $\bar \Q$.
\end{claim}
\begin{proof}
Suppose $\bar p \in \bar \Q^*$. 
Repeatedly strengthen $\bar p$, applying the previous claim for each $\bar s \in {}^{n_{\bar p}}({}^{n_{\bar p}} 2)$
in turn. After finitely many steps we arrive at a condition $\bar q \in D_O$.
\end{proof}
Now let $B = \bigcap_{n\in\omega} O_n$, where each $O_n$ is open dense in $X$. 
By standard arguments, we may find a filter $G$ on $\bar \Q^*$ meeting every $D_{O_n}$ for $n\in \omega$, 
and such that for every $k\in\omega$ there is $\bar p \in G$ with $n_{\bar p} \geq k$ and such that 
$\ran(f_{\bar p(k)})$ has size $2^{n_{\bar p}}$.
For each $k\in\omega$, we obtain a perfect tree 
$T_k = \{ t \setdef (\exists \bar p \in G) (\exists s \in {}^{n_{\bar p(k)}} 2)\; t \subseteq f_{\bar p(k)} (s) \}$;
let $C_k = [T_k]$.
\end{proof}
Note that the previous theorem has the following trivial (but useful!)  corollary; following \cite[p.~274]{blass}, 
one might call this a `polarized version' of Mycielski's Theorem.
\begin{cor}\label{c.myc}
If $k \leq \omega$ and $R$ is a meager $k$-ary relation on ${}^\omega({}^\omega 2)$ 
(i.e.\ $R \subseteq {}^k({}^\omega({}^\omega 2))$, there is for each $i\in k$ a sequence 
$\langle C^i_n \colon n\in \omega\rangle$ perfect subsets of ${}^\omega 2$ such that
$\prod_{i\in k} \prod_{n\in\omega} C^i_n$ is $R$-discrete.
\end{cor}

\subsection{A Galvin-type theorem for iterated Sacks forcing}\label{s.galvin}

Building upon the work in the previous section, we are now ready to give a \emph{topological proof} of 
Galvin's Theorem for iterated Sacks forcing.
We will prove Theorem~\ref{t.galvin} in the following, equivalent form.
It is more easily stated for functions $c^*\colon [\bar p]^2 \to 2$, that is in the context of ordered instead of unordered pairs.

We shall use the analogue of the partition from Definition~\ref{d.delta} in the context of ordered pairs:
\begin{dfn} 
For $\xi <\lambda$ let 
$\Delta^*_\xi = \{ (\bar x,\bar y) \in  ({}^\lambda({}^\omega 2))^2 \setdef \Delta(\bar x,\bar y)=\xi\}$.
\end{dfn}

Striving for generality, we phrase Ramsey theoretic statements for conditions in $\Dt$; 
again, we call to the readers attention that they may replace this set by $\Ds$, if they wish, throughout.
\begin{thm}[Galvin's Theorem for iterated Sacks forcing, 2nd form]\label{t.galvin.second.form}
Let $\bar p\in \Dt[H]$ and $c^*\colon [\bar p]^2 \to 2$ be $\Cantorspace$-universally Baire.
Then there is $\bar q \in \Ds$  with the property that $[\bar q] \subseteq [\bar p]$, 
and for each $\xi \in \dom(H)$, $n \in \omega$, and
$(\bar t_0, \bar t_1) \in \init^n(\bar q)$ such that 
$\bar t_0 \res \sigma = \bar t_1\res\sigma$ and $\bar t_0(\sigma)\neq \bar t_1(\sigma)$,
$c^*$ is constant on
$\big( [\bar q_{\bar t_0}] \times [\bar q_{\bar t_1}]\big) \cap \Delta^*_\xi$.
\end{thm}
This theorem clearly implies Galvin's Theorem for iterated Sacks forcing as stated before, that is, Theorem~\ref{t.galvin}.

We also point out that Theorem \ref{t.galvin.second.form} has a simply corollary for sets (a similar result is proved in 
\cite{miller-nykos} for Borel sets) which we shall use in \S\ref{s.discrete}:
\begin{cor}\label{c.galvin.unary}
Let $\bar p \in \Dt$ and $B \subseteq [\bar p]$ be  $\sigma(\Sigma^1_1)$  (absolutely 
$\Delta^1_2$, $\Cantorspace$-universally Baire, $\textsf{\textup{Y}}^{\apo}$-measurable).
Then there is $\bar q \in \Ds$ such that $[\bar q] \subseteq [\bar p]$ and $[\bar q] \subseteq B$ or 
$[\bar q]\cap B= \emptyset$.
\end{cor}
For this, simply apply Theorem \ref{t.galvin} for 
\[
c(\{ \bar x_0, \bar x_1 \}) = \begin{cases} 1 &\text{ if $\lo(\{ \bar x_0, \bar x_1 \} ) \in B$},\\
0 & \text{ otherwise}
\end{cases}
\]
and restrict to a basic open set.
If we assume $\bar p \in \Ds[\Sigma]$, then we can additionally demand $\bar q \in \Ds[\Sigma]$ 
(but we shall have no use for this).

\medskip

As shall become clear, this solution is optimal and can be arrived at by taking into account obstructions to finding 
homogeneous conditions. 
As our point of departure, we cite the following theorem for continuous colorings (equivalently, clopen partitions) 
due to Geschke, Kojman, Kubi\'s and Schipperus:
\begin{thm}[cf.\ Lemma 37 from \cite{geschke-ea}]\label{t.geschke}
Let $\bar p \in \apo$
and let
$
c\colon \big[ [\bar p] \big]^2 \to 2
$
be a continuous coloring.
Then 
there is $\bar q \in \Ds$, $\bar q\leq\bar p$ such that $c$ is constant on
$\big[ [\bar q] \big]^2$.
\end{thm}
Note that in the language of symmetric maps $c^* \colon [\bar p]^2\to 2$, the the requirement on the coloring becomes 
that $c^* \res \big( [\bar p]^2\setminus \operatorname{diag}([\bar p]^2)\big)$ be continuous.
Clearly, it would be desirable to weaken the requirement of being continuous in the above, 
e.g.\ to being \emph{Borel} or even \emph{Baire measurable}.
The family of colorings $c_\xi$ from \S\ref{s.intro} represent a fundamental obstruction to such a result:
\begin{fct}\label{f.obstruction.1}
Assume momentarily that $\lambda \geq 2$.
Then for any $H$ and $\bar q \in \Dt[H]$ whose support has size at least $2$ and any $\xi \in \supp(\bar q)$, 
$c_{\xi}$ takes both colors on pairs from $[\bar q]^H$.
\end{fct}
This goes towards explaining the role of $\Delta_\xi$.
Moreover, note that $\Delta_0$ is open dense and $\Delta_\xi$ is meager for $\xi>0$. 
Therefore, it is impossible to replace \emph{$\Cantorspace$-universally Baire} by \emph{Baire measurable} 
in Theorem~\ref{t.galvin.second.form}---as a Baire measurable map $c$ can take completely arbitrary values 
on a meager set.

\medskip

Perhaps we can strengthen Theorem~\ref{t.geschke} and show that there exists
$\bar q$ such that when restricted to $\big[[\bar q]\big]^2$, $c$ is constant on each set $\Delta_\xi$ for $\xi<\lambda$,
that is, $c$ becomes a function only of $\Delta(\bar x_0, \bar x_1)$?
An elaboration of the argument from Fact~\ref{f.obstruction.1} uncovers yet another, more subtle obstruction, 
showing that the answer to this question is negative:
\begin{exm}\label{e.obstruction.2}
Let $\bar p \in \apo$; we define a coloring $c$ on $X={}^{\supp(\bar p)}({}^\omega 2)$.
Fix a surjection $g\colon \supp(\bar p) \to \supp(\bar p) \times \omega$.
Also fix functions $\lambda$ and $n$ so that for all $\xi \in \supp(\bar  p)$, $g(\xi)=( \lambda(\xi), n(\xi))$. 

Let $\{\bar x_0, \bar x_1\} \in [X]^2$ be given and let $\xi = \Delta(\bar x_0, \bar x_1)$.
Suppose further that $\bar x_0(\xi)$ comes before $\bar x_1(\xi)$ in the lexicographical ordering. 
Define
\[
c(\{\bar x_0, \bar x_1\}) =
\bar x_1 (\lambda(\xi)) (n(\xi)).
\]
Given $H$ and $\bar q\in \Dt[H]$ stronger than $\bar p$, $c$ obviously induces a coloring of pairs from 
$[\bar q]^{H}$ 
which we also denote by $c$. 
Now suppose that for $\{\bar x, \bar y\} \in \big[[\bar q]\big]^2$,
$c(\{\bar x, \bar y\})$ is simply a function of $\Delta(\bar x, \bar y)$; that is, suppose we can find 
$F \colon \supp(\bar p) \to \{0,1\}$ such that for any $\{\bar x, \bar y\} \in \big[[\bar p]\big]^2$,
$c(\{\bar x, \bar y\})= F( \Delta(\bar x, \bar y))$.

This leads to a contradiction: 
Let $\bar z \in {}^{\supp(\bar p)}({}^\omega 2)$ be given by 
\[
\bar z(\nu)(k) = F( g^{-1} (\nu,k) ).
\]
Find $\bar x \neq \bar z$ such that for each $\xi$,
$\bar x(\xi)$ is the not the right-most branch of $[\bar q]\vert_{\bar x\res \xi}$ (see Definition~\ref{d.top}).
But for any $\xi$, we may chose $\bar y \in [\bar q]$ such that $\Delta(\bar x,\bar y) = \xi$ and $\bar y(\xi)$ comes before 
$\bar x(\xi)$ lexicographically, and thus 
\[
\bar x\big(\lambda(\xi)\big)\big(n(\xi)\big) = c(\{ \bar x, \bar y\}) = F(\xi) = \bar z\big(\lambda(\xi)\big)\big(n(\xi)\big);
\]
contradicting the choice of $\bar x$.
\end{exm}
Taking into account this last obstruction, we arrive at the formulation of Theorem~\ref{t.galvin.second.form}.

\medskip

It will be practical to use not \emph{$\Cantorspace$-universally Baire measurability} but a slightly weaker notion 
of measurability, 
crafted to suit its precise employment in the argument.
\begin{dfn}
Let $\bar p \in \Dt[H]$ and $\xi \in \supp(\bar p)$. We say a set $B \subseteq [\bar p]^H$ is 
$\mathsf{Y}^{\apo}_\xi$-measurable 
if and only if for any $\bar p_i \in \Dt[H_i]$, where $i\in\{0,1\}$, 
such that 
$\bar p_0\res\xi=\bar p_1 \res\xi$, 
$\dom(H_i) = \dom(H)$ 
and $\bar p_i \leq \bar p$ for each $i\in\{0,1\}$ there are 
$\bar q_i \in \Dt[G_i]$ with $\bar q_i \leq \bar p_i$ 
and $\dom(G_i) = \dom(H)$ 
such that $\bar q_0 \res \xi = \bar q_1 \res \xi$ and 
\[
\big( [\bar q_0]^{G_0}\times [\bar q_1 ]^{G_1} \big)\cap \Delta^*_\xi \subseteq B
\]
or
\[
\big ( [\bar q_0]^{G_0}\times [\bar q_1 ]^{G_1} \big ) \cap \Delta^*_\xi \cap B = \emptyset.
\]
Moreover, given $\bar p \in \Dt$ we say a set $B \subseteq [\bar p]^2$ is $\mathsf{Y}^{\apo}$ 
\emph{measurable (on $[\bar p]$)}
just if $B$ is $\mathsf{Y}^{\apo}_\xi$-measurable for every $\xi \in \dom(H)$, for some $H$ such that $\bar p \in \Dt[H]$. 
We say a map $c^* \colon [\bar p]^2 \to X$ for a topological space $X$ is 
$\mathsf{Y}^{\apo}_\xi$-measurable (or $\mathsf{Y}^{\apo}$-measurable) 
just if the pre-image under $c^*$ of every open set in $X$ has said property.
The definition carries over as usual for maps $c\colon \big[[\bar p]\big]^2\to X$ into a topological space $X$.
Note that $c\colon \big[[\bar p]\big]^2\to X$ is $\mathsf{Y}^{\apo}$-measurable 
(resp.\ $\mathsf{Y}^{\apo}_\xi$-measurable) 
\emph{iff} the same holds for the induced symmetric map $c^*$ on $[\bar p]^2$. 
\end{dfn}

We first show $\mathsf{Y}^{\apo}$-measurability is indeed a weaker property than being universally Baire.
This lemma is a crucial part of the proof of our Ramsey theoretic result.
\begin{lem}\label{l.baire=>Y}
Let $\bar p \in \Dt$. Every universally Baire (or just $\Cantorspace$-universally Baire) set 
$B \subseteq [\bar p]^2$ is $\mathsf{Y}^{\apo}$-measurable. Analogously for $B \subseteq \big[[\bar p]\big]^2$.
\end{lem}
Note we are not even assuming $B$ to be symmetric in the first case. 
The proof depends on the `Polarized Mycielski' corollary of Theorem~\ref{t.myc} and on 
a special case of Lemma~\ref{l.Dt.hom}\eqref{l.second} (and its proof).
\begin{proof}
Let $\bar p \in \Dt[H]$ and $\bar p_i \in \Dt[H_i]$ with $\dom(H_0) = \dom(H)$ for each $i\in \{0,1\}$
be given.
We may assume that 
\[
\big(\forall (\bar x_0, \bar x_1) \in [\bar p_0]\times [\bar p_1]\big)\; \bar p_0(\sigma)\neq \bar p_1(\sigma).
\]
It suffices to treat the case $B \subseteq [\bar p]^2$, for we may replace
$B \subseteq \big[[\bar p]\big]^2$ by a symmetric subset of $[\bar p]^2$.

Write $D = \dom(H)$, $\alpha = \otp D$, write $\pi$ for the order preserving map $\pi\colon D\to\alpha$ 
and write $\Q ={}^\alpha \Sacks$.
For the moment, fix $i\in\{0,1\}$.
Let
\[
h_i\colon [\bar p_0] \to {}^\alpha(\Cantorspace) = [1_{\Q}]
\]
be the ``obvious'' homeomorphism; for concreteness, it is defined as follows:

Given any perfect tree $T\subseteq {}^{<\omega}2$, write
$h_T \colon T \to  {}^{<\omega}2$ for the unique bijection which sends  $(T)^*_n$ to ${}^{<\omega}2$
and preserves the lexicographic ordering; this gives us a homeomorphism
\[
h^*_T \colon  [T]  \to\Cantorspace.
\]
For each $\bar x \in [\bar p_i]$ we define $h_i(\bar x) \in  {}^\alpha(\Cantorspace)$ by recursion on $\xi < \alpha$, 
letting
\[
h_i(\bar x) (\xi) = h^{*}_T\big(\bar x(\xi')\big),
\] 
with $\xi' = \pi^{-1}(\xi)$ and $T = (H_i)_{\xi'}(\bar x\res \xi')$.
Writing $\bar p$ for $\bar p_0 \res \sigma$, which equals $\bar p_1 \res \sigma$, 
note that $h_0 \res [\bar p] = h_1 \res [\bar p]$.

We therefore obtain a homeomorphism
\[
\big({}^\alpha(\Cantorspace) \times {}^\alpha(\Cantorspace)\big)\cap\Delta^*_{\pi(\sigma)} 
\stackrel{h_0 \times h_1}{\cong}
 \big([\bar p_0]\times [\bar p_1]\big) \cap\Delta^*_\sigma
\]
The space on the right-hand side is equipped with the subspace topology. 
Writing 
\[
\gamma := \alpha - \pi(\sigma), 
\]
we clearly we have another homeomorphism, denoted by $g$,
\begin{equation}\label{e.g}
{}^{\pi(\sigma)}(\Cantorspace) \times {}^{\gamma}(\Cantorspace) \times {}^{\gamma}(\Cantorspace)
\stackrel{g}{\longrightarrow} \big({}^\alpha(\Cantorspace) \times {}^\alpha(\Cantorspace)\big)\cap\Delta^*_{\pi(\sigma)}
\end{equation}
the precise form of which we leave to the reader to guess correctly.
Consider 
\[
B' = {\big(g \circ (h_0 \times h_1)\big)^{-1}}'' B,
\]
a Baire measurable set. Either $B'$ or its complement must be non-meager;
assume $B'$ is non-meager (in the other case, the argument is almost identical).
Since $D$ is countable, the left-hand side in \eqref{e.g} is homeomorphic to ${}^\omega\Cantorspace$,
and we may apply the polarized generalization of Mycielski's Theorem (that is, Corollary~\ref{c.myc})
and obtain a product of perfect sets,
\[
\bar C:=\prod_{\xi < \pi(\sigma)} C_\xi 
\times \prod_{\zeta < \gamma} C^0_\zeta 
\times \prod_{\eta < \gamma} C^1_\eta
\subseteq B'
\]
Applying $g$ we obtain two conditions, $\bar q'_0, \bar q'_1 \in {}^\alpha \Sacks$ such that 
$\bar q'_0 \res \pi(\alpha) = \bar q'_1 \res \pi(\alpha)$ so that $g$ restricts to a homeomorphism on $\bar C$,
\[
\bar C = \prod_{\xi < \pi(\sigma)} C_\xi 
\times \prod_{\zeta < \gamma} C^0_\zeta 
\times \prod_{\eta < \gamma} C^1_\eta
\stackrel{g \res \bar C} {\longrightarrow}
\big([\bar q'_0] \times [\bar q'_1]\big)\cap\Delta^*_\sigma,
\]
and moreover,
\begin{equation}\label{e.subset.B}
\big({h_0}''[\bar q'_0] \times {h_1}''[\bar q'_1]\big)\cap\Delta^*_\sigma \subseteq B.
\end{equation}
We may treat $\bar q_i :=1_{\Q}$ and $\bar q'_i$ (for $i\in\{0,1\}$) as very simple iterated Sacks conditions in which all names are standard names
for ground model perfect sets; we even have $\bar q_i, \bar q'_i \in \Dt$.
Then by Lemma~\ref{l.Dt.hom}\eqref{l.second}, for each $i\in\{0,1\}$ we can find $\bar p'_i \in \Dt$
such that 
\[
{h_i} ''  [\bar q'_i] = [\bar p_i']
\]
and by the proof of said Lemma~\ref{l.Dt.hom}\eqref{l.second} also $\bar p'_0 \res \sigma = \bar p'_1 \res\sigma$.
From \eqref{e.subset.B} we conclude 
\[
\big([\bar p'_0]\times [\bar p'_1]\big) \cap\Delta^*_\sigma \subseteq B
\]
as required.
In case $B'$ is meager, clearly the same argument goes through with $B$ replaced by its 
complement, $[\bar p]^2 \setminus B$.
\end{proof}

\medskip

Now
we can prove the main theorem of this section. 
\begin{proof}[Proof of Theorem~\ref{t.galvin.second.form}]
We use our standard fusion argument:
Suppose $\bar p \in \Dt$ and $c^*\colon [\bar p ]^2 \to \{0,1\}$ is 
$\textsf{\textup{Y}}^{\apo}$-measurable (it suffices to consider such $c^*$ by Lemma~\ref{l.baire=>Y}).
We may also assume $\bar p\in \Ds[\Sigma]$.

As in the previous lemma, let $\bar p_0 =\bar p$ and build a sequence 
$\bar p_0 \geq_1 \bar p_1 \geq_2 \bar p_2 \geq_3 \hdots$ ($\Sigma$ remains fixed throughout the argument).
Assuming we have constructed $\bar p_{n}$, we specify how to obtain 
$\bar p_{n+1} \leq_{n+1} \bar p_{n}$. 
Observe the following:
\begin{claim}\label{c.thin.above}
Let $\bar s^0, \bar s^1 \in \init^{n+1}(\bar p_n)$ such that $\bar s^0 \neq \bar s^1$,
and $\bar q \in \apo$ such that $\bar q \leq_{n+1} \bar p_n$ be given, and
let $\xi$ be least such that 
$\bar s^0(\xi) \neq \bar s^1(\xi)$.
There is $\bar q^* \leq_{n+1} \bar q$ such that
$c^*$ is constant on $\big([\bar q^*_{\bar s^0}] \times[\bar q^*_{\bar s^1}]\big)\cap\Delta^*_\xi$.
\end{claim}
\begin{proof}[Proof of Claim]
Using $\mathsf{Y}^{\apo}$-measurability of $c^*$ find 
$\bar r^0, \bar r^1 \in \apo$ such that $\bar r^i \leq \bar q_{\bar s}$ for each $i\in\{0,1\}$, 
$\bar r^0 \res \xi = \bar r^1 \res\xi$, and $c^*$ is 
constant 
on $\big([\bar r^0] \times[\bar r^1]\big)\cap\Delta^*_\xi$. 
Now ``thin out'' $\bar q$ above 
$\bar s^0$ and $\bar s^1$, to obtain $\bar q^*$ as in the claim:
\begin{subclaim}
There is $\bar q^* \in \apo$ such that 
$\bar q^* \leq_{n+1} \bar q$ and $\bar q^*_{\bar s^i} = \bar r^i$ for each $i\in\{0,1\}$.
\end{subclaim}
\noindent
\textit{Proof of Subclaim.}
This is technical but straightforward. 
We provide details nevertheless; 
we cannot hinder a trustful reader from skipping the proof of the subclaim.

Write $\bar s$ for $\bar s^0\res\xi$ (which equals $\bar s^1\res\xi$)
and $\bar r$ for $\bar r^0\res\xi$ (which equals $\bar r^1\res\xi$).

Construct $\bar q^*$ by recursion on $\sigma <\lambda$.
Suppose we have already constructed $\bar q^* \res \sigma$.
For $\sigma < \xi$ find $\bar q^*(\sigma)$ such that
\[
(\bar q^*\res\sigma)_{\bar s\res\sigma}\forces \bar q^*(\sigma)_{\bar s(\sigma)} = \bar r(\sigma)
\]
and
\[
(\bar q^*\res\sigma)_{\bar t\res\sigma}\forces \bar q^*(\sigma)_{\bar t(\sigma)} = \bar q(\sigma)_{\bar t(\sigma)}
\]
whenever $\bar t \in \init^{n+1}(\bar q)$ and $\bar t\res\sigma \neq \bar s\res\sigma$.
This is possible because 
$(\bar q^*\res\sigma)_{\bar t\res\sigma}$ and $(\bar q^*\res\sigma)_{\bar s\res\sigma}$
 are incompatible.

We now treat the case $\sigma = \xi$: 
Find $\bar q^*(\xi)$ such that for each $i \in \{0,1\}$
\[
(\bar q^*\res\xi)_{\bar s}\forces \bar q^*(\xi)_{\bar s^i(\xi)} = \bar r^i(\xi)
\]
which is possible because $s^0(\xi) \neq s^1(\xi)$, and such that
\[
(\bar q^*\res\xi)_{\bar t\res\xi}\forces \bar q^*(\xi)_{\bar t(\xi)} = \bar q(\xi)_{\bar t(\xi)}
\]
whenever $\bar t \in \init^{n+1}(\bar q)$ and $\bar t\res\xi \neq \bar s\res\xi$.
The latter is possible for the same reason as in the case $\sigma<\xi$.

Finally, for $\sigma>\xi$, find $\bar q^*(\sigma)$ such that for each $i \in \{0,1\}$
\[
(\bar q^*\res\sigma)_{\bar s^i\res\sigma}\forces \bar q^*(\sigma)_{\bar s^i(\sigma)} = \bar r^i(\sigma)
\]
which is possible because $s^0\res\sigma \neq s^1\res\sigma$, and such that
\[
(\bar q^*\res\sigma)_{\bar t\res\sigma}\forces \bar q^*(\sigma)_{\bar t(\sigma)} = \bar q(\sigma)_{\bar t(\sigma)}
\]
whenever $\bar t \in \init^{n+1}(\bar q)$ and $\bar t\res\sigma \notin \{\bar s^0\res\sigma, \bar s^1\res\sigma\}$.
The latter is possible for the same reason as previously.
\renewcommand{\qedsymbol}{{\tiny Subclaim and Claim.} $\Box$}
\end{proof}
Now build a finite $\leq_{n+1}$-descending sequence of conditions $\bar q^0$, $\bar q^1$, \dots below $\bar p_n$, 
applying the claim for each pair $\bar s^0, \bar s^1 \in \init_{n+1}(\bar p_n)$.
Let $\bar p_{n+1}$ be the last element of the sequence.
We let $\bar q$ be the greatest lower bound of the sequence $\langle \bar p_n \setdef n\in \omega\rangle$; 
by construction, 
$c^*$ is continuous on each subspace $\Delta^*_\xi$ in the strong sense required in 
Theorem~\ref{t.galvin.second.form}.
\renewcommand{\qedsymbol}{{\tiny Theorem~\ref{t.galvin.second.form}.} $\Box$}
\end{proof}

\section{An absoluteness result}\label{s.abs}

We shall need the following proposition for our proof of Theorem~\ref{t.main},
our application of Galvin's Theorem in the realm definable of maximal discrete sets.
\begin{prop}\label{p.elementarity}
Let $\bar s$ be $\apo$-generic over $\Vee$ and suppose $\lambda \geq \omega_1$ 
(recalling that throughout this article, $\apo$ denotes a countable support iteration of Sacks forcing and 
$\lambda$ is its length).
Then $\Vee[\bar s\res \omega_1] \prec_{\Sigma^1_3} \Vee[\bar s]$.
\end{prop}
For the proof, we shall need the following.
\begin{lem}\label{l.equiv}
As usual, let $\apo$ be a iteration of Sacks forcing of any length. 
Suppose $\bar p \in \Dt[H]$ and $F_y \colon [\bar p]_{H} \to \Bairespace$ 
and $\Phi$ is $\Pi^1_1$. 
Then 
\begin{equation}\label{e.p.forces}
\bar p \forces_{\apo} (\forall z \in \Bairespace)\; \Phi(F_y(\bar s^{\dot G}),z)
\end{equation}
if if and only if
\begin{equation}\label{e.p.below}
(\forall \bar p' \in \Dt[H']\text{ s.t. }\bar p'\leq \bar p)(\forall F_z \colon [\bar p']_{H'} \to \Bairespace)
(\exists \bar u \in [\bar p']_{H'})\; 
\Phi\big(\bar F_y(\bar u), F_z(\bar u)\big)
\end{equation}
where,
letting $\proj \colon [\bar p'] \to {}^{\supp(\bar p)}(\Cantorspace) $ denote the map 
$\bar x \mapsto \bar x \res \supp(\bar p)$,
$\bar F_y$ denotes the unique map such that $\bar F_y = F_y \circ \proj$.
\end{lem}
We shall use each of the directions of the above equivalence for a Sacks iteration of different length.
\begin{proof}
We show the contrapositives.
Supposing first \eqref{e.p.below} fails, find $\bar p'$ and $F_z$ witnessing its failure, that is, 
such that $\bar p' \in \Dt$, $\bar p' \leq\bar p$ and
\begin{equation*}\label{e.forall.u}
(\forall \bar u \in [\bar p'])\; \neg \Phi(\bar F_y(\bar u),F_z(\bar u)).
\end{equation*}
Since this formula is $\Pi^1_2$ it is forced to hold in $\Vee[\bar s]$, and since 
$\bar p' \forces \bar s^{\dot G} \in [\bar p']$, we conclude 
$\bar p' \forces \neg \Phi\big(\bar F_y(\bar s^{\dot G}), F_z(\bar s^{\dot G})\big)$,
and since $\bar p' \leq \bar p$, \eqref{e.p.forces} must fail.

For the other direction, suppose \eqref{e.p.forces} fails, whence there is
$\bar p' \leq \bar p$ such that
$\bar p' \forces (\exists z \in \Bairespace)\; \neg\Phi(\bar F_y(\bar s^{\dot G}),z)$, 
where $\bar F_y$ is as in the statement of the lemma.
By Lemma~\ref{l.continuous.reading.gen} (and replacing $\bar p'$ by a stronger condition, if needed) 
we can assume that $\bar p' \in \Dt$ and find $F_z\colon [\bar p'] \to \Bairespace$ 
such that $\bar p' \forces  \neg\Phi\big(F_y(\bar s^{\dot G}),F_z(\bar s^{\dot G})\big)$.
Finally, 
letting $B=\big\{\bar u \in [\bar p'] \setdef \Phi\big(F_y(\bar u),F_z(\bar u)\big)\big\}$ and
applying Corollary~\ref{c.galvin.unary} we can assume (replacing $\bar p'$ by a stronger condition if needed)
that $[\bar p'] \cap B = \emptyset$, that is, 
$(\forall \bar u \in [\bar p']) \neg \Phi\big(F_y(\bar u),F_z(\bar u)\big)$.
Thus \eqref{e.p.below} fails. 
\end{proof}

Given Lemma~\ref{l.equiv} we can prove Proposition~\ref{p.elementarity}.
\begin{proof}
Let $\Psi$ be the formula
\[
(\exists y \in \Bairespace)(\forall z \in \Bairespace)\; \Phi(y,z)
\] 
where $\Phi$ is $\Sigma^1_1$, and suppose $\Vee[\bar s] \vDash \Psi$.
We can assume any parameter of $\Psi$ is an element of $\Vee$ and 
suppress it.
Find $H$ and $\bar p \in \Dt[H]$ and $F_y \colon [\bar p]_{H} \to \Bairespace$
such that
\begin{equation*}
\bar p \forces_{\apo} (\forall z \in \Bairespace)\; \Phi(F_y(\bar s),z).
\end{equation*}
By the previous lemma, \eqref{e.p.below} holds.

Now let $\pi\colon \supp(\bar p) \to \alpha$, where $\alpha \in\omega_1$, denote the transitive collapse.
Write
\[
\Q := \apo_{\omega_1}
\]
and find 
\[
\bar q \in \Q \cap \Dt 
\]
and a homeomorphism $h$ as in Lemma~\ref{l.Dt.hom},
\[
h\colon [\bar q] \to [\bar p].
\]
Moreover, let
\[
F'_y = F_y \circ h.
\]
We shall now show the following:
\begin{claim}\label{c.below}
It must hold that
\begin{equation}\label{e.q.below}
(\forall \bar q' \in \Dt \cap \Q \text{ s.t. }\bar q'\leq \bar q)(\forall F'_z \colon [\bar q'] \to \Bairespace)
(\exists \bar u \in [\bar q']_{H'})\; 
\Phi\big(\bar F'_y(\bar u), F'_z(\bar u)\big)
\end{equation}
where in analogy to Lemma~\ref{l.equiv}, 
recalling that $\alpha = \supp(\bar q)$ and letting 
$\proj'\colon [\bar q'] \to {}^\alpha(\Cantorspace)$ denote the projection $\bar x \mapsto \bar x \res \alpha$,
we define
$\bar F'_y := F'_y \circ \proj'$.

In other words, \eqref{e.p.below} must hold with $\apo$, $\bar p$, and $F_y$ replaced 
by $\Q$, $\bar q$, and $F'_y$ respectively. 
\end{claim}
\begin{proof}[Proof of claim]
Suppose towards a contradiction \eqref{e.q.below} fails and find $\bar q' \in \Dt$ and $F'_z$ witnessing this failure, 
that is, such that
$\bar q' \leq \bar q$ and
\begin{equation}\label{e.q.witness}
(\forall \bar u \in [\bar q'])\; 
\neg\Phi\big(\bar F'_y(\bar u), F'_z(\bar u)\big).
\end{equation}
As in Lemma~\ref{l.Dt.hom}, find  $\bar p' \in \apo \cap \Dt$ such that $\bar p' \leq \bar p$ and a homeomorphism 
\[
h'\colon [\bar q'] \to [\bar p']
\]
In fact, we can require that
$\big[\bar p' \res \sup \supp(\bar p)\big] \subseteq [\bar p]$ and that 
\[
h' = h \circ \proj'.
\]
Let 
\[
F_z := F'_z \circ (h')^{-1}.
\]
We leave it to the reader to verify that $\bar F'_y = \bar F_y \circ h'$.
Then \eqref{e.q.witness} entails 
\[
(\forall \bar u \in [\bar p'])\; 
\neg\Phi\big(\bar F_y(\bar u), F_z(\bar u)\big),
\]
contradicting \eqref{e.p.below}. Having reached a contradiction, \eqref{e.q.below} must hold.
\renewcommand{\qedsymbol}{{\tiny Claim.} $\Box$}
\end{proof}
By \eqref{e.q.below} and Lemma \ref{l.equiv}, we have 
\[
\bar q \forces_{\Q} (\forall z \in \Bairespace)\; \Phi(F'_y(\bar s^{\dot G}),z),
\]
and thus we have found $\bar q \in \Q$ such that $\bar q \forces \Psi$.
Since $\Q$ is weakly homogeneous, any formula is which is 
forced by any condition $\bar q\in\Q$ is forced by the trivial condition $1_{\Q}$.
Therefore $1_{\Q} \forces \Psi$, proving the proposition.
\renewcommand{\qedsymbol}{{\tiny Proposition~\ref{p.elementarity}.} $\Box$}
\end{proof}

\section{Definable discrete sets in the iterated Sacks extension}\label{s.discrete}

We shall now prove Theorem~\ref{t.main} in slightly more general form:
\begin{thm}\label{t.main.gen}
Let $a \in {}^\omega \omega$, let $\mathcal R$ be a $\Sigma^1_1[a]$ binary relation on $\omega^\omega$, 
and let $\bar s$ be generic sequence of reals for an iteration of Sacks forcing over $\eL[a]$.
Then there is a $\Sigma^1_2[a]$-formula 
which defines a maximal $\mathcal R$-discrete set in both $\eL[a]$ and $\eL[a][\bar s]$. 
(In fact, this set  is $\Delta^1_2[a]$, witnessed by the same pair of formulas in both models).
\end{thm}

First, we draw a corollary from the proposition of the previous section:

\begin{cor}\label{c.omega_1}
It is enough to prove Theorem~\ref{t.main.gen} for the iteration of Sacks forcing of length
$\omega_1$, that is, in the case $\lambda=\omega_1$.
\end{cor}
\begin{proof}
Suppose $\bar s$ is the generic sequence for an iterated Sacks extension of $\eL[a]$ of arbitrary length.
If Theorem~\ref{t.main.gen} holds for $\lambda=\omega_1$, we can fix a $\Sigma^1_2(a)$ formula
$\Psi(x,a)$ such 
\[
\eL[a,\bar s\res\omega_1]\vDash\text{``}\{x\in{}^\omega\omega \setdef \Psi(x,a)\}\text{ is maximal $\mathcal R$-discrete.''}
\]
Since the formula on the right can be expressed by a $\Pi^1_3$ statement, it persist in $\eL[a,\bar s]$ 
by Proposition~\ref{p.elementarity}.
\end{proof}

Lastly, we will use the following rather natural terminology: 
\begin{dfn}\label{d.galvin.witness}
 Let $\bar p \in \Dt$, a continuous function $F\colon [\bar p] \to {}^\omega \omega$ and a binary relation 
$\mathcal R$ on ${}^\omega \omega$ be given.
\begin{enumerate}[label=(\alph*)]
 \item We say $\bar q \leq p$ is 
a \emph{Galvin witness (for $(F,\mathcal R)$) iff} 
the conclusion of Theorem~\ref{t.galvin.second.form} holds when $c^*$ is the
characteristic function on $[\bar p]^2$ of $\mathcal R$, that is, \emph{iff}
$\bar q \in \Ds$ and 
$[\bar q]\subseteq [\bar p]$ and for each $\sigma \in \supp(\bar q)$, $n \in \omega$, and 
$(\bar t_0, \bar t_1) \in \init^n(\bar q)$ such that $\bar t_0 \res \sigma = \bar t_1\res\sigma$ and 
$\bar t_0(\sigma)\neq \bar t_1(\sigma)$,
\begin{equation}\label{e.complete.triple}
\big( F''[\bar q_{\bar t_0}] \times F''[\bar q_{\bar t_1}]\big) \cap \Delta^*_\sigma \subseteq \mathcal R
\end{equation}
or
\begin{equation}\label{e.discrete.triple}
\big( F''[\bar q_{\bar t_0}] \times F''[\bar q_{\bar t_1}]\big) \cap \Delta^*_\sigma \cap \mathcal R = \emptyset.
\end{equation}
\item 
We say that $(\sigma,\bar t_0, \bar t_1)$ is \emph{complete in $\bar q$} 
(or more precisely, \emph{$(F,\mathcal R)$-complete in $\bar q$}) if
$\sigma \in \supp(\bar q)$, 
$(\bar t_0, \bar t_1) \in \init^n(\bar q)$,  for some $n \in \omega$, 
$\bar t_0 \res \sigma = \bar t_1 \res \sigma$
 and 
$\bar t_0(\sigma)\neq \bar t_1(\sigma)$,
and
\eqref{e.complete.triple} holds.
\item 
Similarly, we call a triple $(\sigma,\bar t_0, \bar t_1)$ as in the previous item \emph{discrete in $\bar q$} 
(or more precisely, \emph{$(F,\mathcal R)$-discrete in $\bar q$}) if
\eqref{e.discrete.triple} holds
instead of
\eqref{e.complete.triple}.
\end{enumerate}
\end{dfn}
It is a trivial consequence of Theorem~\ref{t.galvin} that such a Galvin witness can always be found. 
The reader may find it instructive to note at this point that ``$\bar q$ is a Galvin witness'' is a 
$\Pi^1_2$ statement in an appropriate code for $(\bar q,\Sigma)$.
We find it probably that this can be lowered to $\Delta^1_2$ (or better) by thorough analysis, 
but have been able to avoid 
the need for this elaboration in the arguments below.

\medskip

Let $\bar p$, $F$, and $\mathcal R$ be fixed as in Definition~\ref{d.galvin.witness} for now. 
We shall also make use of the following observations:
\begin{fct}\label{f.obs}
 Suppose $\bar q \leq \bar p$ is a Galvin witness.
\begin{enumerate}
\item 
One of the following holds:
\begin{enumerate}[label=(\Alph*)]
\item $F''[\bar q]$ is $\mathcal R$-discrete,
\item\label{i.galvin.witness.complete}
There is a complete triple $(\sigma,\bar t_0, \bar t_1)$ in $\bar q$.
\end{enumerate}
\item\label{f.i.discrete} Suppose $(\sigma,\bar t_0, \bar t_1)$ is complete in $\bar q$
and there is no complete triple in $\bar q$ with first component $\sigma' <\sigma$.
Then for any $\sigma' < \sigma$, 
\[
(F \times F)''\big([\bar q]^2 \cap \Delta^*_{\sigma'}\big) \cap \mathcal R = \emptyset.
\]
\end{enumerate}
\end{fct}
Suppose now that $(\sigma,\bar t_0, \bar t_1)$ is a complete triple in $\bar q$
whose first component is minimal among complete triples.

\medskip

We are now ready to prove Theorem~\ref{t.main.gen} (and thus Theorem~\ref{t.main}).
 
\begin{proof}[Proof of Theorem~\ref{t.main.gen}]
By Corollary~\ref{c.omega_1} we may assume $\lambda = \omega_1$.
The proof relativizes easily to the parameter $a$, so we suppress it below. 
Therefore, let us suppose $\mathcal R$ is a $\Sigma^1_1$ binary relation on $\Bairespace$ given by
\[
u \mathbin{\mathcal R} v \iff (\exists w)\; u \mathbin{\mathcal{R}^w} v
\]
where $u \mathbin{\mathcal{R}^w} v$ is a $\Sigma^0_\omega$ relation of triples $(u,v,w)$.

Fix an enumeration $\langle (\bar p_\xi,  \Sigma_\xi, F_\xi) \setdef \xi <\omega_1\rangle$ of all triples 
$(\bar p, \Sigma,  F)$ where $\bar p \in \Ds[\Sigma]$,
$\supp(\bar p) \in\lambda$, 
and $F \colon [\bar p] \to {}^\omega \omega$ is continuous, with the triples enumerated in
 increasing order with respect to $\leq_{\eL}$.

We will construct a sequence of trees $\langle T_\xi \setdef \xi<\omega_1\rangle$ such that
$[T_\xi] \subseteq [\bar p_\xi]^{\Sigma_\xi}$.
For each $\xi \leq \omega_1$, let $\dot{\mathcal A}_\xi$ be a name such that
\[
1_{\apo}\forces\dot{\mathcal A}_\xi = \bigcup_{\nu <\xi}  {F_\nu}''[\check T_\nu].
\]
and write $\dot {\mathcal A}$ for $\dot{\mathcal A}_{\omega_1}$
We will take care to ensure that:
\begin{enumerate}[label=(\roman*)]
\item\label{i.A.discete} It is forced that $\dot{\mathcal A}_{\omega_1}$ is maximal $\mathcal R$-discrete.
\item\label{i.A.definable} The set $\{T_\xi \setdef \xi<\omega_1\}$ is a $\Sigma^1_2$, and hence $\mathcal A$ is forced to be $\Sigma^1_2$
as well.
\end{enumerate}

The construction of $\langle T_\xi \setdef \xi<\omega_1\rangle$ is by recursion.
Suppose we have already constructed $\langle T_\nu \setdef \nu < \xi\rangle$.

Let me give a short, rough sketch how we shall find $T_\xi$. 
Firstly, we find $\bar p' \leq \bar p_\xi$ which decides
\begin{equation}\label{sketch.e.forces.new}
\text{``}\{ F_\xi(\bar s_{\dot G}) \}\cup\dot{\mathcal A}_\xi\text{ is $\mathcal R$-discrete'',}
\end{equation}
setting $T_\xi = \emptyset$ if $\bar p'$ forces the negation of \eqref{sketch.e.forces.new}.
If $\bar p'$ forces \eqref{sketch.e.forces.new} to hold,
we will be able to find $T_\xi$ such that 
\begin{itemize}
\item ${F_\xi}''[T_\xi]$ is $\mathcal R$-discrete, 
and 
\item there is $\bar p'' \leq \bar p_\xi$ which forces that
$\{F_\xi(\bar s_{\dot G})\} \cup {F_\xi}''[\check T_\xi]$ is not $\mathcal R$-discrete.
\end{itemize}
It is not hard to see that this ensures \ref{i.A.discete};
to also ensure \ref{i.A.definable} we have to choose $T_\xi$ as above very carefully, namely as follows.

Find 
the $\leq_{\eL}$-least 
quadruple $(T,\bar p',  \Sigma', \vec{w})$
 such that 
$\bar p' \leq \bar p_\xi$, $\bar p' \in \Ds[\Sigma']$,
and
one of the following conditions \eqref{i.first} or \eqref{i.second} holds:

\begin{enumerate}[label=(\Roman*), ref=\Roman*]
\item\label{i.first} $T=\emptyset$ and $\vec{w} = (Y,W)$ where 
$Y\colon [\bar p'] \to \bigcup_{\nu <\xi}  [T_\nu]$ and 
$W\colon [\bar p'] \to \Bairespace$
are continuous, and
\begin{equation}\label{e.discrete}
(\forall \bar x \in [\bar p'])  \; F_\xi(\bar x) \mathbin{\mathcal{R}^{W(\bar x)}} (F_\nu \circ Y)(\bar x),
\end{equation}
or
\item\label{i.second} 
$(\forall \bar x \in [p']) (\forall \nu < \xi) (\forall \bar y \in [T_\nu])  \;  F_\xi(\bar x) \nR F_\nu(\bar y)$, 
and moreover, $\vec w = (\bar p'',\Sigma'',Y,W)$, where $\bar p'' \in \Ds[\Sigma'']$ and
all of the following hold:
\begin{enumerate}
\item $\bar p''\leq \bar p'$,
\item \label{ii.discrete} ${F_\xi}''[T]$ is $\mathcal R$-discrete,
\item $Y\colon [\bar p''] \to [T]$ and $W\colon \colon [\bar p''] \to \Bairespace$ are continuous,  and
\item \label{ii.complete} 
$[\bar p'']\subseteq [T]$ or
\begin{equation}\label{e.complete}
\big(\forall \bar x \in[\bar p'']\big)  \; F_\xi(\bar x) \mathbin{\mathcal{R}^{W(\bar x)}} (F_\xi \circ Y)(\bar x);
\end{equation}
\end{enumerate}
\end{enumerate}
We will show in a moment such a quadruple does indeed exists\footnote{Moreover
the reader may verify that
in fact,
the possibilities \eqref{i.first},  \eqref{ii.discrete}, and \eqref{ii.complete} are mutually exclusive.
Since we take a minimal quadruple, we do not need this.}.
Having chosen $(T,\bar p',  \Sigma', \vec{w})$ to be the $\leq_{\eL}$-least one, 
define
\[
T_\xi = T.
\]
We must show that this is well-defined:
\begin{claim}\label{c.witness}
There is a quadruple $(T,\bar p',  \Sigma', \vec{w})$ satisfying \eqref{i.first} or \eqref{i.second} as above.
\end{claim}
\begin{proof}[Proof of claim]
We argue by cases.

First assume $\bar p_\xi \not\forces$``$\dot{\mathcal A}_\xi\cup\{F_\xi(\bar s_{\dot G})\}$ is $\mathcal R$-discrete''.
Then we can find $\bar q \leq \bar p_\xi$, a $\apo$-names $\dot{\bar y}$ and $\dot w$, and $\nu < \xi$ such that
$\bar q \forces$``$\dot{\bar y} \in [T_\nu]$ and $F_\nu(\dot{\bar y}) \mathbin{\mathcal R^{\dot w}} F_\xi(\bar s_{\dot G})$''.
By Lemma~\ref{l.continuous.reading.gen} and  we may find $\bar q' \leq \bar q$ 
such that $\bar q' \in \Ds$ 
and 
continuous functions $Y \colon [\bar q'] \to [T_\nu]$
and $W \colon [\bar q'] \to \Bairespace$ such that
that $\bar q' \forces$``$\dot{\bar y} = Y(\bar s_{\dot G})$ and $\dot w = W(\bar s_{\dot G})$'', whence
\begin{equation}\label{e.bar.q'.forces.R}
\bar q' \forces (F_\nu\circ Y)(\bar s_{\dot G}) \mathbin{\mathcal R^{W(\bar s_{\dot G})}} F_\xi(\bar s_{\dot G}).
\end{equation}
Let
\[
B =\{\bar x \in [\bar q'] \setdef F_\xi(\bar x) \mathbin{\mathcal R^{W(\bar x)}} (F_\nu \circ Y)(\bar x)\}
\]
 and apply Corollary~\ref{c.galvin.unary} to find $\bar p' \leq \bar q'$ such that $\bar p' \in \Ds$ and $[\bar p'] \subseteq B$
 or $[\bar p'] \cap B =\emptyset$. 
By \eqref{e.bar.q'.forces.R} and by upwards absoluteness of $\Sigma^1_1$ formulas, 
it must be the case that $[\bar p'] \subseteq B$; in other words, \eqref{e.discrete} holds, from which 
it is clear that 
$(T,\bar p',  \Sigma', \vec{w})$ as required in \eqref{ii.complete} exists.

We proceed to the second case, 
namely that $\bar p_\xi \forces$``$\dot{\mathcal A}_\xi\cup\{F_\xi(\bar s_{\dot G})\}$ is $\mathcal R$-discrete''.
Letting
\[
B =\{\bar x \in [\bar p_\xi] \setdef \big(\forall \nu < \xi\big)\big(\forall \bar y \in [T_\nu]\big) 
 F_\xi(\bar x) \nR F_\nu(\bar y)\}
\]
we can argue much like in the previous case to obtain $\bar q' \leq \bar p$ such that $\bar q' \in \Dt$ and
$\bar q' \subseteq B$.
Finally, find $\bar p' \leq \bar q$ such that $\bar p' \in \Ds$ and such that $\bar p'$ is a 
Galvin witness for $F_\xi$. 

We must distinguish to sub-cases;
The first sub-case to consider is when ${F_\xi}''[\bar p']$ is $\mathcal R$-discrete.
In this case, find $T$ and $\vec{w}$ as required in \eqref{i.second} as follows:
Simply let $\bar p'' =\bar p'$ and find $T$ such that $[T]=[\bar p']$.
Since the first alternative in \eqref{ii.complete} obtains, it is clear that \eqref{i.second} is satisfied
with $\vec w = (p'', Y, W)$ where $Y$ and $W$ can be chosen arbitrarily.

The second and remaining sub-case is that ${F_\xi}''[\bar p']$ is not $\mathcal R$-discrete.
Let $\sigma$ be least such that there exist a complete triple in $\bar p'$ whose first component is
$\sigma$, and let $(\sigma,\bar t_0, \bar t_1)$ be such a triple.
Note that then by minimality of $\sigma$ and since $\bar p'$ is a Galvin witness,
\begin{equation}\label{e.choice.of.sigma}
(\forall \sigma' < \sigma)\; \big({F_\xi}''[\bar p']\big)^2 \cap \Delta^*_{\sigma'} \cap \mathcal R = \emptyset
\end{equation}
Next, we find $Y$ as in \eqref{i.second}.
\begin{subclaim}\label{sc.Y}
For any $\bar x \in [\bar p'_{\bar t_0}]_\sigma$ we can find $Y(\bar x) \in [\bar p'_{\bar t_1}]$ such that
$Y(\bar x)\res \sigma = \bar x$, and so that $\bar x \mapsto Y(\bar x)$ is a continuous function.
\end{subclaim}

\begin{proof}[Proof of Subclaim]
We give the details for the convenience of the reader. Define 
\[
Y(\bar x) = \bar y
\]
to be the unique left-most minimal element of $[\bar p'_{\bar t_1}]$ 
such that $\bar y\res \sigma = \bar x$.
Or, to be pedantically explicit,  
with $H$ such that $\bar p'_{\bar t_1} \in \Dt[H]$ define
$\bar y \in [\bar p'_{\bar t_1}]$ recursion: 
$\bar y\res \sigma = \bar x$
and for each $\sigma' \in\dom(H)\setminus \sigma$,
$\bar y(\sigma')$ is  the left-most branch of $H_{\sigma'}(\bar y \res \sigma')$,
and $\bar y(\sigma')$ is the left-most branch of ${}^{<\omega}2$ for $\sigma' \notin\dom(H)\cup\sigma$.
\renewcommand{\qedsymbol}{{\tiny Subclaim \ref{sc.Y}.} $\Box$}
\end{proof}
As $\bar t_0(\sigma)$ and $\bar t_1(\sigma)$ are $\subseteq$-incomparable,
$\big(\bar x, Y(\bar x\res\sigma)\big) \in \Delta^*_\sigma$
for any $\bar x \in [p'_{\bar t_0}]$ and so
$\bar x  \mathbin{\mathcal R}  Y(\bar x \res \sigma)$.
By $\mathbf{\Pi}^1_1$-absoluteness,
we infer
\[
\bar p'_{\bar t_0} \forces 
 F_\xi(\bar s_{\dot G}) \mathbin{\mathcal R} (F_\xi\circ Y)(\bar s_{\dot G}\res\sigma).
\]
Argue as previously (in the first case, and using Lemma~\ref{l.continuous.reading.gen} ) 
to find $\bar p'' \leq \bar p'_{\bar t_0}$ and $W$ as required in \eqref{ii.complete}. 
Finally, find $T$ such that $[T]=[\ran(Y)]$.

\begin{subclaim}\label{sc.discrete}
The set $ {F_\xi}''[T]$ is $\mathcal R$-discrete.
\end{subclaim}

\begin{proof}[Proof of Subclaim.]
Let $y_0, y_1 \in {F_\xi}''[T]$ such that $y_0 \neq y_1$. Since $[T] = Y''[\bar p'_{\bar t_0}]_\sigma$
we can find $\bar y_i \in [\bar p]_\sigma$ such that $y_i = (F_\xi\circ Y)(\bar y_i)$ for each $i\in\{0,1\}$,
and clearly $\bar y_0 \neq \bar y_1$. 
As $Y(\bar y_i) \res \sigma = \bar y_i$, $\big(Y(\bar y_0), Y(\bar y_1)\big) \in \Delta^*_{\sigma'}$ 
for some $\sigma'<\sigma$,
and so $y_0 \nR y_1$ by \eqref{e.choice.of.sigma}.
\renewcommand{\qedsymbol}{{\tiny Subclaim \ref{sc.discrete}.} $\Box$}
\end{proof}
It is now clear how to find $(T,\bar p',  \Sigma', \vec{w})$ as required in \eqref{i.second}.
\renewcommand{\qedsymbol}{{\tiny Claim~\ref{c.witness}.} $\Box$}
\end{proof}

This finishes our construction of $\langle T_\xi \setdef \xi<\omega_1\rangle$.
Supposing $\bar s$ is $(\apo,\eL)$-generic and working in $\eL[\bar s]$, we can now define 
$\mathcal A$ as follows:
\[
\mathcal A := \{x\in\Bairespace \setdef (\exists \xi<\omega_1)\; x \in {F_\xi}''[T_\xi]\},
\]
that is, $\mathcal A = \dot{\mathcal A}^G$.

\begin{claim}\label{c.discrete}
The set $\mathcal A$ is $\mathcal R$-discrete.
\end{claim}
\begin{proof}[Proof of claim]
Let distinct elements $x, y \in \mathcal A$ be given and find $\xi, \nu < \omega_1$ 
and $\bar x$, $\bar y$ such that
such that 
$x = F_\xi(\bar x)$, $y = F_\nu(\bar y)$, and
$\bar x \in [T_{\xi}]$ and $\bar y \in [T_{\nu}]$.

We distinguish two cases:
Firstly, if $\xi \neq \nu$, we can assume $\nu < \xi$. 
Note that since $T_{\xi} \neq \emptyset$, \eqref{i.second} must obtain, and so 
by \eqref{e.discrete} and because $\mathbf{\Pi}^1_1$-formulas are absolute, we infer $x \nR y$.

Suppose otherwise that $\xi = \nu$, in which case $\{\bar x, \bar y\}\subseteq [T_\xi]$.
In this case $x \nR y$ is immediate by \eqref{ii.discrete} and 
because $\mathbf{\Pi}^1_1$-formulas are absolute.
\renewcommand{\qedsymbol}{{\tiny Claim~\ref{c.discrete}.} $\Box$}
\end{proof}

\begin{claim}\label{c.maximal}
The set $\mathcal A$ is maximal among $\mathcal R$-discrete sets.
\end{claim}
\begin{proof}[Proof of claim]
Suppose towards a contradiction that $\bar p \forces$``$\dot{\mathcal A}$ is not maximal $\mathcal R$-discrete''.
We can assume (by Lemma~\ref{l.continuous.reading.gen}) that 
$\bar p \in \Ds$ and
for some continuous $F \colon [\bar p] \to \Bairespace$,
\begin{equation}\label{e.intruder}
\bar p \forces \big\{F(\bar s_{\dot G})\big\}\cup\dot{\mathcal A}\text{ is $\mathcal R$-discrete.}
\end{equation}
Moreover by Theorem~\ref{t.galvin.second.form} we can assume that $\bar p$ is a Galvin witness for $F$.

We can find $\xi<\omega_1$ such that $(\bar p, F) = (\bar p_\xi, F_\xi)$.
Let $(\bar p', \bar \Sigma', T, \vec w)$ be $\leq_{\eL}$-least such that
\eqref{i.first} or \eqref{i.second} holds (we have seen in Claim~\ref{c.witness} that such a quadruple exists).
But \eqref{i.first} is impossible since \eqref{e.discrete} contradicts \eqref{e.intruder},
and so \eqref{i.second} obtains.

We distinguish two cases, namely whether 
the first or second alternative in \eqref{ii.complete} 
obtains.
In the first case, by construction
$[\bar p'']\subseteq [T_\xi]$.
But then, recalling $F=F_\xi$, and again using $\mathbf{\Pi}^1_1$-absoluteness,
\[
\bar p'' \forces F(\bar s_{\dot G}) \in {F_\xi}''[T_\xi]\subseteq \dot{\mathcal A} 
\]
contradicting the choice of $\bar p$ and $F$.

In the second and remaining case,  
\eqref{e.complete} holds,
a formula $\mathbf{\Pi}^1_1$ and hence absolute.
Since $\bar p'' \forces \bar s_{\dot G} \in [p_{\bar t_1}]$, 
we have
\[
\bar p'' \forces F(\bar s_{\dot G}) \mathbin{\mathcal R} (F_\xi \circ Y) (\bar s_{\dot G}) \in \dot{\mathcal A}.
\]
Again, this contradicts the choice of $\bar p$ and $F$.
\renewcommand{\qedsymbol}{{\tiny Claim~\ref{c.maximal}.} $\Box$}
\end{proof}

\begin{claim}\label{c.definable}
The set $\mathcal A$ is $\Sigma^1_2$.
\end{claim}
\begin{proof}[Proof of claim]
We show that
\begin{align}\label{e.def.A}
x \in \mathcal A \iff (\exists \beta < \omega_1)
\eL_\beta\vDash\text{``}\mathsf{ZF}^-  \land \; (\exists \xi < \omega_1)\; x\in {F_\xi}''[T^{\eL_\beta}_\xi]
\text{''}
\end{align}
where $T_\xi^{\eL_\beta}$ denotes the $\xi$th element of the sequence $\bar T^{\eL_\beta}$ 
obtained as the outcome of interpreting the definition of 
$\bar T = \langle T_\nu \setdef \nu <\omega_1\rangle$ in $\eL_\beta$.
From \eqref{e.def.A} it follows by the usual methods that $\mathcal A$ is $\Sigma^1_2$, 
since the statement on the right
is known to be equivalent to a $\Sigma^1_2$ formula with free variable $x$.

To prove \eqref{e.def.A} it suffices to show that 
\begin{equation}\label{e.define.T}
(\xi, T) \in \bar T \iff (\exists \beta \in \omega_1 \setminus \xi)
\eL_\beta\vDash\text{``}\mathsf{ZF}^- \text{'' and } (\xi, T) \in \bar T^{\eL_\beta}.
\end{equation}
As is routine to verify, \eqref{i.first} and \eqref{i.second} are each expressible by 
a formula which is at most $\mathbf{\Pi}^1_1$ in an appropriate code for $(\bar p', \Sigma', \vec{w})$,
and thus \eqref{i.first} and \eqref{i.second} are each equivalent to a $\Sigma_1$-formula in any
transitive $\in$-model of $\mathsf{ZF}^-$ containing $(\bar p', \Sigma', \vec{w})$ as an element, 
and absolute for such models.
Moreover $\langle \bar p_\xi, F_\xi\setdef\xi < \omega_1\rangle$ is definable in a way which is absolute 
for initial segments of $\eL$.
From this, standard argument show that $\bar T$ is definable by a $\Sigma_1$ recursion, 
hence itself $\Sigma_1$, and also \eqref{e.define.T} holds.
\renewcommand{\qedsymbol}{{\tiny Claim~\ref{c.definable}.} $\Box$}
\end{proof}
This finishes the proof of the theorem.
\renewcommand{\qedsymbol}{{\tiny Theorem~\ref{t.main.gen}.} $\Box$}
\end{proof}

\section{A co-analytic maximal orthogonal family in the iterated Sacks extension}\label{s.mof}

In \cite{fischer-tornquist}, the authors show that the space $P(X)$ of Borel probability measures on a 
perfect effective Polish space $X$ is itself an effective Polish space, and that the orthogonality relation is Borel.
The same holds for $P_c(X)$, the Borel subset  of atomless measures in $P(X)$.
It follows immediately from Theorem~\ref{t.main} that in the iterated Sacks extension of length 
$\omega_2$ there is a $\Delta^1_2$ mof in $P_c(X)$.
As any measure with an atom is non-orthogonal to a Dirac measure we immediately obtain 
Theorem~\ref{t.mof}. by  quoting the following result from \cite{schrittesser-tornquist}:
\begin{thm}
If $X$ is a perfect effective Polish space and 
there is a $\Delta^1_2$ (equivalently, a $\Sigma^1_2$) mof in $P_c(X)$ then there is a $\Pi^1_1$ mof in $P_c(X)$.
\end{thm}

We round off the results from the present paper as well as from \cite{schrittesser-tornquist} 
with an observation due to B.\ Miller. We see this as evidence that there is no substantially simpler way of finding 
mofs in forcing extension with new reals, e.g.\ one cannot find a mof in $\eL$ that is indestructible by forcing.
\begin{thm}\label{t.miller}
Let $M$ be an inner model of $\Vee$ such that $\Power(\omega)^M \neq \Power(\omega)$ and let 
$X$ be a perfect Polish space with code in $M$. 
Then there is $\mu \in P(X)$ such that $(\forall \nu \in M) \; \mu \perp \nu$.
In particular, if $\mathcal A \in M$ and $M\vDash \mathcal A$ is a maximal orthogonal families of measures in 
$P(X)$ then $\mathcal A$ is not maximal in $\Vee$.
\end{thm}
\begin{proof}
For the purpose of this proof, identify $\Cantorspace$ with a subset of $X$.
Consider the map
$F\colon \Cantorspace^2 \to \Cantorspace$ given by
\[
F(x,y) (n) = \begin{cases}
x(\frac{n}{2}) &\text{ if $n$ is even, }\\
y(\frac{n-1}{2}) & \text{ if $n$ is odd.}
\end{cases}
\]
for $n\in \omega$.
This gives rise to a perfect family of perfect sets
$\langle P_x \setdef x \in \Cantorspace \rangle$, where $P_x =\{F(x,y) \setdef y \in \Cantorspace \}$.

For any $\nu \in P(X)^M$, since the measure algebra has the ccc there is a countable 
$H\subseteq \Cantorspace$  in $M$ such that 
\begin{equation*}
 (\forall x \in \Cantorspace) \; x \in H \iff \nu(P_x)> 0.
\end{equation*}
This statement is $\Pi^1_1(H)$ and hence absolute between $M$ and $\Vee$.
In particular, if $y \in \Cantorspace^{\Vee} \setminus M$ and $\nu \in M$ it must be the case that $\nu(P_y)  = 0$.

Now suppose $y \in \Cantorspace^{\Vee} \setminus M$. 
Find by the product measure construction $\mu \in P(X)$ such that $\mu(P_y) = 1$. 
For any $\nu \in M$ we have $\nu \perp \mu$.
\end{proof}

\section{Questions}\label{s.questions}

I shall conclude with the following two questions which are, to my knowledge, open.

\begin{enumerate}[label=\arabic*.]
\item Is there a model in which there is a $\Pi^1_1$ mof but \emph{no} $\Pi^1_1$ mad family?
\item Is there a $\Pi^1_1$ (or equivalently, a $\Sigma^1_2$) mof in the Laver extension?
\end{enumerate}

\bibliography{mof-large-continuum}{}
\bibliographystyle{amsplain}

\end{document}